\documentclass[10pt]{article}
\usepackage[dvips]{graphics}
\usepackage{epsfig,graphicx}
\usepackage{amsfonts}
\usepackage{amsmath}

\input xy
\xyoption{all}

\begin{document}
\newcommand{\lar}{\longrightarrow}
\newcommand{\kF}{\mathcal F}
\newcommand{\kQ}{\mathcal Q}
\newcommand{\kG}{\mathcal G}
\newcommand{\kL}{\mathcal L}
\newcommand{\kB}{\mathcal B}
\newcommand{\kM}{\mathcal M}
\newcommand{\kO}{\mathcal O}
\newcommand{\kP}{\mathcal P}
\newcommand{\kI}{\mathcal I}
\newcommand{\kC}{\mathcal C}
\newcommand{\kJ}{\mathcal J}
\newcommand{\bJ}{\mathbf J}
\newcommand{\bE}{\mathbf E}
\newcommand{\bF}{\mathbf F}
\newcommand{\kE}{\mathcal E}
\newcommand{\kT}{\mathcal T}
\newcommand{\kN}{\mathcal N}
\newcommand{\kS}{\mathcal S}
\newcommand{\kk}{\mathbf k}
\newcommand{\kA}{\mathcal A}
\newcommand{\PP}{\mathbf P}
\newcommand{\ZZ}{\mathbb Z}
\newcommand{\RR}{R}
\newcommand{\balpha}{\mathbf \alpha}
\newcommand{\bbetha}{\mathbf \beta}
\newcommand{\FF}{\mathbf F}
\newcommand{\tens}{\otimes}
\newcommand{\Hom}{\mathop{\rm Hom}}
\newcommand{\Coh}{\mathop{\rm Coh}}
\renewcommand{\mod}{\mathop{\rm mod}}
\newcommand{\TF}{\mathop{\rm TF}}
\newcommand{\T}{\mathop{\rm T}}
\newcommand{\End}{\mathop{\rm End}}
\newcommand{\Ext}{\mathop{\rm Ext}\nolimits^1}
\newcommand{\TC}{\mathop{\rm TC}}

\newtheorem{theorem}{Theorem}[section]
\newtheorem{corollary}[theorem]{Corollary}
\newtheorem{remark}[theorem]{Remark}
\newtheorem{example}[theorem]{Example}
\newtheorem{lemma}[theorem]{Lemma}
\newtheorem{proposition}[theorem]{Proposition}
\newtheorem{statement}[theorem]{Statement}
\newtheorem{definition}[theorem]{Definition}
\newtheorem{conjecture}[theorem]{Conjecture}
\newtheorem{problem}[theorem]{Problem}
\newtheorem{question}[theorem]{Question 1}

\title
{
Indecomposables of the derived categories of certain associative algebras
 }

\author
{
Igor Burban 
\footnote {Kaiserslautern University,  Kyiv University and Institute of Mathematics }  
\and Yurij Drozd 
\footnote{Kyiv Taras Shevchenko University and Kaiserslautern University}
}
\date{}

\maketitle

\begin{abstract}
In this article we describe indecomposable objects of the derived categories of a branch class of  associative algebras. To this class belong such known 
classes of algebras as gentle algebras, skew-gentle algebras and certain 
degenerations of tubular algebras.   
\end{abstract}

\setcounter{footnote}{0}
\renewcommand{\thefootnote}{}
 \footnote{
The both authors were partially supported CRDG Grant UM2-2094. The first 
author was also partially supported by DFG Schwerpunkt ``Globale Methoden 
in der komplexen Geometrie''.}

\section{\bf Introduction} 
In this article we elaborate a method of description of  indecomposable 
objects of derived categories of representations of associative algebras proposed in \cite{Jalgebra}. We consider the derived categories of gentle and skew-gentle algebras treated also in 
\cite{Skowr}, \cite{Ringel1},
\cite{laPena}, by using the Happel's functor from the derived category of finite-dimensional modules over a finite-dimensional $\kk$-algebra $A$ of finite homological dimension to the stable category over the repetitive algebra 
$\hat{A}$. This class of algebras was also treated in  \cite{Bekkert1} and 
\cite{Bekkert2} by using other methods. 
 An advantage of our approach is that it also works for algebras of 
infinite homological dimension and gives a description of the derived category $D^{-}(A-\mod)$ of   \emph{bounded from the right} complexes.  

It is well-known that canonical tubular algebras of type $(2,2,2,2, \lambda)$

\begin{figure}[ht]
\hspace{3.9cm}
\includegraphics[height=3cm,width=4.5cm,angle=0]{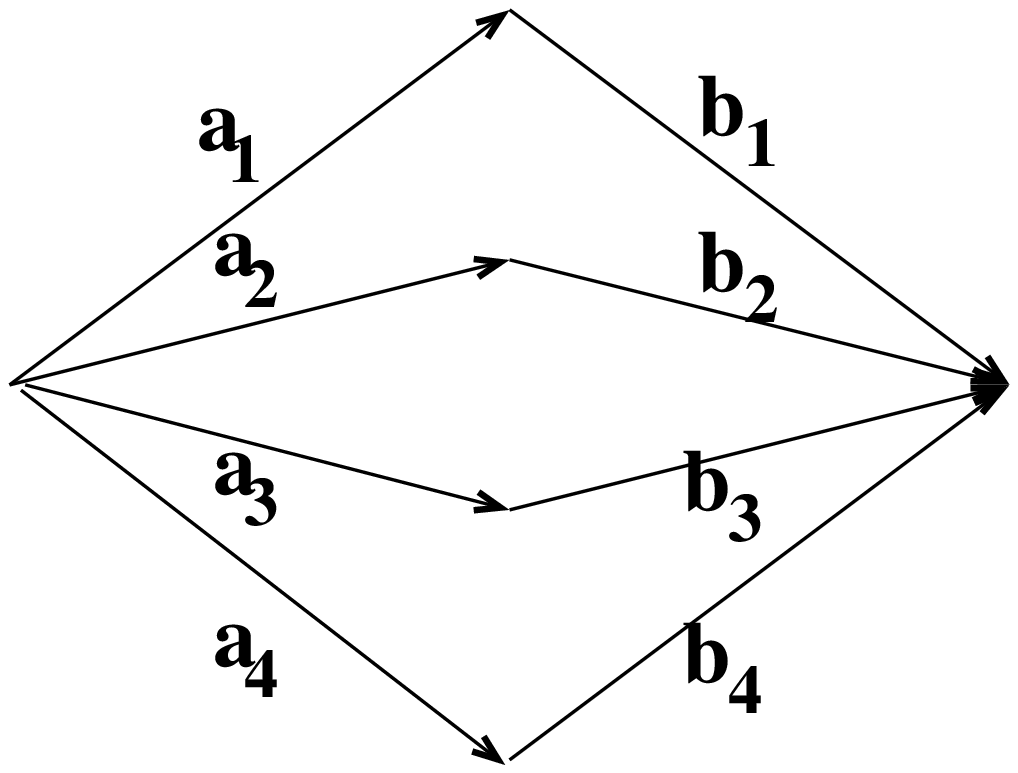}
\end{figure}

\noindent
where relations are 
$$
 b_1 a_1 + b_2 a_2 + b_3 a_3  = 0, 
$$
$$
 b_2 a_2 + b_3 a_3 + \lambda b_4 a_4  = 0 
$$
and  $\lambda \neq 0,1$ are derived-tame of polynomial growth. 
They arises naturally in connection with weighted projective lines of tubular type \cite{Lenzing}. A natural question is: what happens if a family of tubular algebras specializes to a forbidden value of parameter $\lambda = 0$? It turns out that the derived category of the degenerated tubular algebra is still tame but  now it has exponential growth.

\section{\bf Category of triples} 

Let ${\mathcal A}$ and ${\mathcal B}$ be two abelian categories, 
${\bf F}: {\mathcal A}\longrightarrow {\mathcal B}$ a right exact functor. It induces a functor 
between corresponding derived categories 
${\bf D^{-} \bf F}: D^{-}({\mathcal A})\longrightarrow D^{-}({\mathcal B})$. The description of  the fibers of this functor is often equivalent to certain matrix problem.

In the following we are going to work in the following situation. 
Let $A$ be a semi-perfect associative $\kk$-algebra 
(not necessarily finite dimensional),
$A\subset \tilde{A}$ be an embedding such that $r = rad(A) = rad(\tilde{A})$.
Let $I\subset A$ be a  two-sided $\tilde{A}$-ideal  containing $r$. It means 
that 
 $r\subseteq I = I\tilde{A} = \tilde{A} I$, thus
 $A/I$ and $\tilde{A}/I$ are semi-simple algebras.

\medskip
\noindent
\begin{example}
Consider the following embedding of algebras:

\begin{figure}[ht]
\hspace{3.7cm}
\includegraphics[height=3.5cm,width=5cm,angle=0]{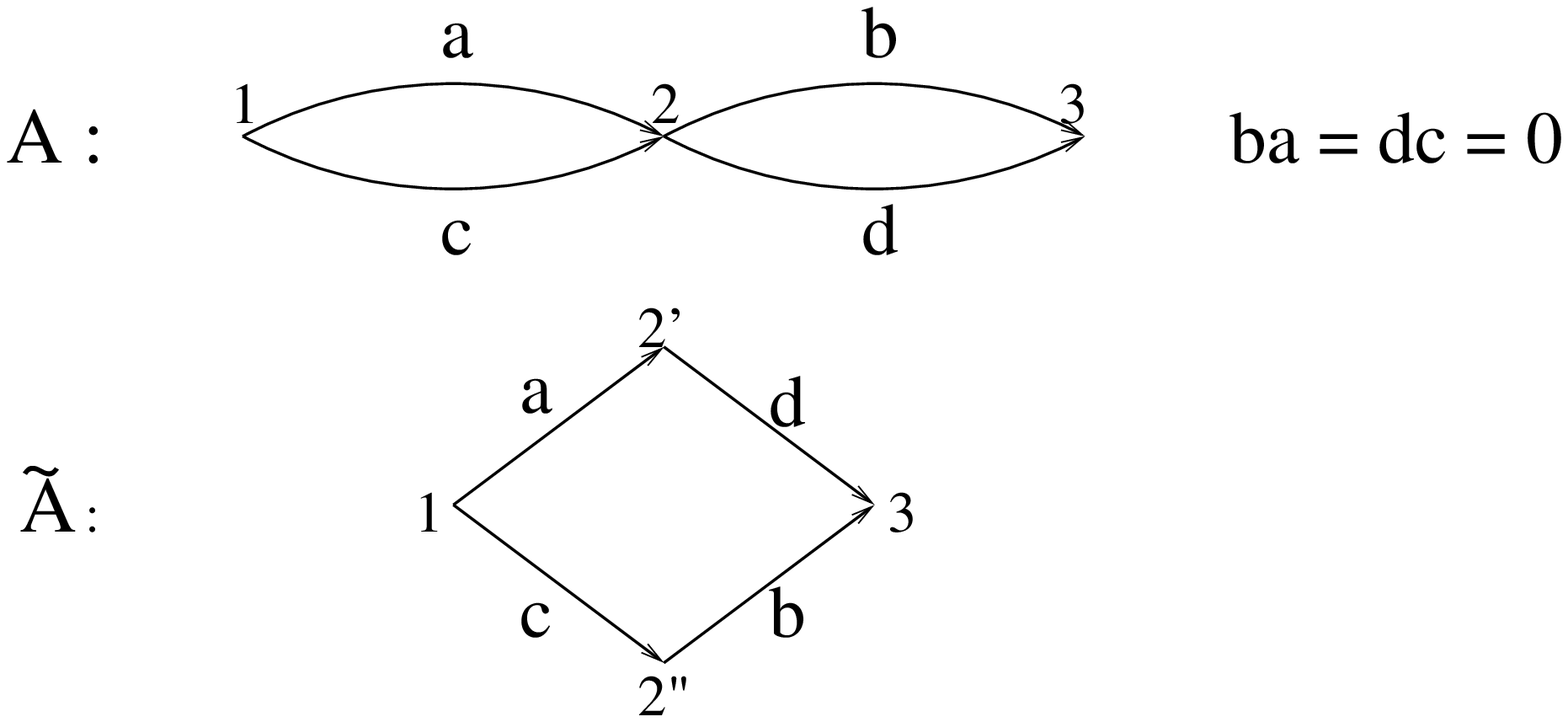}
\end{figure}

\noindent
Then the two-sided ideal $I$ in $A$ and $\tilde{A}$, which is generated by $e_1$ and $e_3$ satisfies 
the properties above. 
\end{example}

\begin{definition}
Consider the following category of triples of 
complexes $\TC_{A}$ 

\begin{enumerate}
\item Objects are  triples 
$(\widetilde {\mathcal P}_{\bullet}, {\mathcal M}_{\bullet}, i)$, where

$\tilde{\mathcal P}_{\bullet} \in D^{-}(\tilde{A}-\mod)$, 

${\mathcal M}_{\bullet} \in D^{-}(A/I-\mod)$, 

$i : {\mathcal M}_{\bullet} \lar
\tilde{A}/I\tens_{\tilde{A}} \tilde{\mathcal P}_{\bullet}$ 
 a morphism in
 $ D^{-}(A/I-\mod) $, such that 

$\tilde{i} :
\tilde{A}/I\tens_{A/I} {\mathcal M}_{\bullet} \lar
\tilde{A}/I\tens_{\tilde{A}}\tilde{\kP}_{\bullet}$
is an isomorphism in $D^{-}(\tilde{A}/I-\mod).$ 

\item Morphisms
 $                              
(\tilde{\mathcal P_{\bullet}}_{1}, {\mathcal M_{\bullet}}_{1}, i_{1}) \lar
(\tilde{\mathcal P_{\bullet}}_{2}, {\mathcal M_{\bullet}}_{2}, i_{2})
 $
are  pairs $(\Phi, \varphi)$, 
%\hspace{3cm}
$$\tilde{\mathcal P_{\bullet}}_{1} \stackrel{\Phi}\lar 
\tilde{\mathcal P_{\bullet}}_{2},\quad
 {\mathcal M_{\bullet}}_{1} \stackrel{\varphi}\lar
 {\mathcal M_{\bullet}}_{2},
$$ such that 

\begin{tabular}{p{2.3cm}c}
 &
\xymatrix
{ 
\tilde{A}/I\tens_{\tilde{A}} \tilde{\mathcal P_{\bullet}}_{1} 
\ar[d]^{
\Phi\tens \mbox{\rm id}} &
{\mathcal M_{\bullet}}_{1} \ar[l]_{\qquad i_1} \ar[d]^{\varphi} \\               \tilde{A}/I\tens_{\tilde{A}} \tilde{\mathcal P_{\bullet}}_{2}  &
{\mathcal M_{\bullet}}_{2} \ar[l]^{\qquad i_2} 
}
\end{tabular}

is commutative. 
\end{enumerate}
\end{definition}

\begin{remark}
If an algebra $A$ has  infinite homological dimension, then we are forced to deal
with the derived category of right  bounded  complexes (in order
to define  the left derived functor of the tensor product). In case $A$ has  
finite homological dimension we can suppose that all complexes above are {\it bounded}
from  both sides.
\end{remark}

\begin{theorem}[\cite{Jalgebra}]
The functor 
$$
D^{-}(A-\mod) \stackrel{\mathbf F}\lar {\rm TC}_{A}
$$
$
\kP_{\bullet} \lar
(\tilde{A}\tens_{A} {\mathcal P}_{\bullet} ,
A/I\tens_{A}{\mathcal P}_{\bullet},
i : A/I\tens_{A} {\mathcal P}_{\bullet}   \lar \tilde{A}/I\tens_{A}{\mathcal P}_{\bullet}) 
$
has the following properties:

\begin{enumerate}
\item ${\mathbf F}$ is dense (i.e., every triple 
$(\tilde{\kP}_{\bullet}, \kM_{\bullet}, i)$ is 
isomorphic to some
${\mathbf F}(\kP_{\bullet})$). 
\item $\FF(\kP_{\bullet}) \cong \FF(\kQ_{\bullet}) \Longleftrightarrow 
\kP_{\bullet} \cong \kQ_{\bullet}$.
\item ${\mathbf F}(\kP_{\bullet})$ is indecomposable if and only if so  is $\kP_{\bullet}$
(note that this property is an easy formal consequence of the previous two properties).
\item ${\mathbf F}$ is full.
\end{enumerate} 
\end{theorem}

The main point to be clarified is: having a triple 
$\kT = (\tilde{\kP}_{\bullet}, \kM_{\bullet}, i)$ how can we reconstruct $\kP_{\bullet}$?
The exact sequence 
$$0\lar I\tilde{\kP}_{\bullet} \lar \tilde{\kP}_{\bullet} \lar
\tilde{A}/I\tens_{\tilde{A}} \tilde{\kP}_{\bullet} \lar 0$$ of complexes
in $A-\mod$ gives a distinguished triangle 
$$
I\tilde{\kP}_{\bullet} \lar \tilde{\kP}_{\bullet} \lar
\tilde{A}/I\tens_{\tilde{A}} \tilde{\kP}_{\bullet}
\lar I\tilde{\kP}_{\bullet}[-1]
$$ 
in $D^{-}(A-\mod)$.
The properties of triangulated categories imply that there is a
 morphism of triangles

\begin{tabular}{p{1.4cm}c}
&
\xymatrix
{
I\tilde{\mathcal P}_{\bullet} \ar[r] &
\tilde{\mathcal P}_{\bullet} \ar[r] &
\tilde{A}/I\tens_{\tilde{A}} \tilde{\mathcal P}_{\bullet}  \ar[r]&
I\tilde{\mathcal P}_{\bullet}[-1] \\
I\tilde{\mathcal P}_{\bullet} \ar[r] \ar[u]^{id} &
{\mathcal P}_{\bullet} \ar[r] \ar[u]_{\Phi} &
{\mathcal M}_{\bullet} \ar[r] \ar[u]_{i} &
I\tilde{\mathcal P}_{\bullet}[-1],\ar[u]^{id}\\
}
\end{tabular}

\noindent
where $\kP_{\bullet} = cone(\kM_{\bullet} \lar I \tilde{\kP_{\bullet}}[-1])[1]$. Set ${\bf G}(\kT) = \kP_{\bullet}.$ Taking a cone
is not a functorial operation. It gives an intuitive explanation why
the functor ${\bf F}$ is not an equivalence.
The properties of triangulated categories  immediately imply that
the constructed map (not a functor!)
$$
{\bf G} : Ob({\rm TC}_{A}) \lar Ob\bigl(D^{-}(A-\mod)\bigr)
$$
sends isomorphic objects into isomorphic ones and 
${\bf GF}(\kP_{\bullet}) \cong \kP_{\bullet}$.  For more details see
\cite{Jalgebra}.

\section{\bf Derived categories of gentle algebras}
It was observed that the representation theory of gentle (or, more general,
string) algebras is closely related to the representation of quivers
of type $\tilde{A}_n$.  We shall sketch a proof of the result of
Z.~Pogorzaly, A.~Skowronski \cite{Skowr} that any gentle algebra is derived-tame and propose a new way  of reduction  to a matrix problem.

Let $A$ be the  path algebra of the quiver

\begin{figure}[ht]
\hspace{3.7cm}
\includegraphics[height=1.5cm,width=5cm,angle=0]{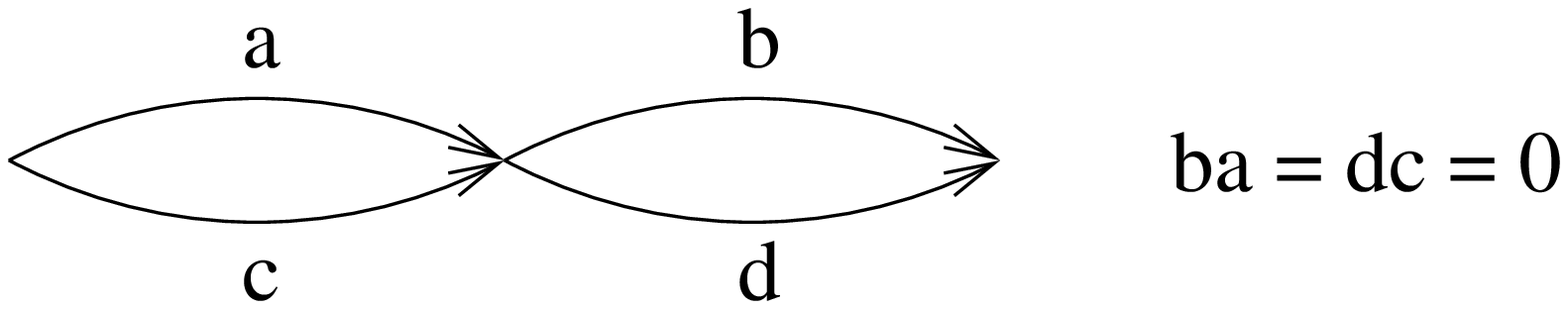}
\end{figure}

Then we can embed it into the path algebra $\tilde{A}$ of the quiver

\begin{figure}[ht]
\hspace{4.7cm}
\includegraphics[height=2.1cm,width=2.3cm,angle=0]{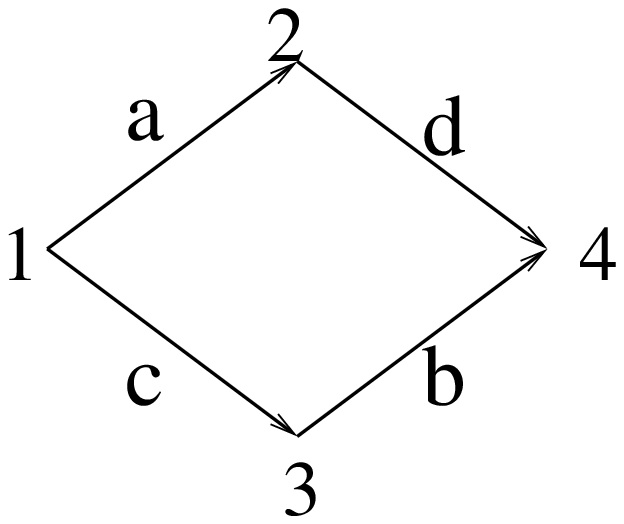}
\end{figure}

In this case set $I = (a,b,c,d,e_{1},e_{4})$. So $A/I = \kk$ and 
$\tilde{A}/I = \kk\times \kk$, $A/I \longrightarrow \tilde{A}/I$ is a diagonal embedding.

As we have seen in the previous section, a complex $\kP_{\bullet}$ of the
derived category $D^{-}(A-\mod)$ is defined by some triple
$(\tilde{\kP}_{\bullet},\kM_{\bullet},i)$. Since $A/I-\mod$ can be identified with
the category of $\kk$-vector spaces, the map 
$i: \kM_{\bullet} \lar \tilde{\kP}_{\bullet}/I\tilde{\kP}_{\bullet}$ is given by a collection
of linear maps
$H_{k}(i) : H_{k}(\kM_{\bullet}) \lar H_{k}(\tilde{\kP}_{\bullet}/I\tilde{\kP}_{\bullet})$.
The map $H_{k}(i)$ is a $\kk$-linear map of a $\kk$-module into a 
$\kk\times \kk$-module.  Hence it is given by two matrices $H_{k}(i)(1)$ and 
$H_{k}(i)(2)$. From the non degenerated condition of the category of
 triples it follows that both of these matrices are square and non degenerated.

The algebra $\tilde{A}$ has a homological dimension $1$.  By the theorem of Dold (see \cite{Dold}), an indecomposable complex
from $D^{-}(\tilde{A}-\mod)$ has a form
$$
\dots \lar 0 \lar \underbrace{M}_{i} \lar 0 \lar \dots, 
$$ 
where $M$ is an indecomposable $\tilde{A}$-module.

The next question is: which transformations can we perform with the matrices
$H_{k}(i)(1)$ and $H_{k}(i)(2)$? 

We can do, simultaneously, any elementary transformation of columns of matrices
$H_{k}(i)(1)$ and $H_{k}(i)(2)$. From the definition of the category of triples
it follows that row transformations are induced by the morphisms in
$D^{-}(\tilde{A}-\mod)$.

If $\tilde{P}_{i}$ denotes the indecomposable projective $\tilde{A}$-module
corresponding to the vertex $i$, then
$$
\tilde{A}/I \otimes_{\tilde{A}} \tilde{P}_{1} =
\tilde{A}/I \otimes_{\tilde{A}} \tilde{P}_{4} = 0,
$$ 
$$
\tilde{A}/I \otimes_{\tilde{A}} \tilde{P}_{2} = \kk(2),\quad
\tilde{A}/I \otimes_{\tilde{A}} \tilde{P}_{3} = \kk(3).
$$

\noindent
Consider the continuous series of representations of the quiver
$\tilde{A} = \tilde{A}_4$:

\begin{figure}[ht]
\hspace{4cm}
\includegraphics[height=2.8cm,width=2.8cm,angle=0]{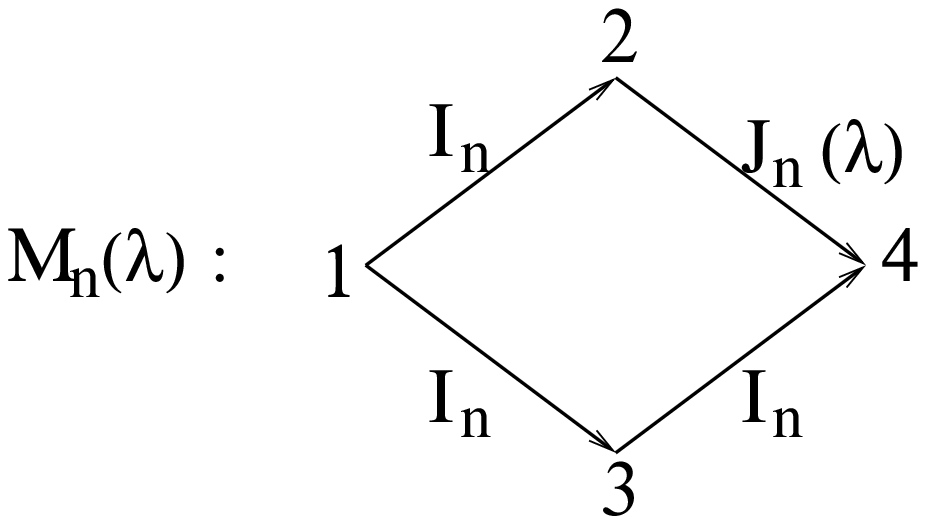}
\end{figure} 

\noindent
$M_{n}(\lambda)$ has a projective resolution:
$$
0 \longrightarrow \tilde{P}_{4}^{n} \longrightarrow \tilde{P}_{1}^{n} \longrightarrow M_{n}(\lambda) \longrightarrow 0.
$$
Hence, in the derived category $D^{-}(\tilde{A}/I)$ holds:
$\tilde{A}/I \otimes_{\tilde{A}} M_{n}(\lambda) = 0$. So, 
$\tilde{A}/I\otimes_{\tilde A}$ {\it kills} the continuous series 
of representations of $\tilde{A}$. We only have to consider the discrete series of representations.

Recall some basic facts of the theory of representations of tame hereditary algebras 
(see \cite{Ringel} for more details). The Auslander-Reiten quiver
has the following structure:

\begin{figure}[ht]
\hspace{3cm}
\includegraphics[height=3.1cm,width=6cm,angle=0]{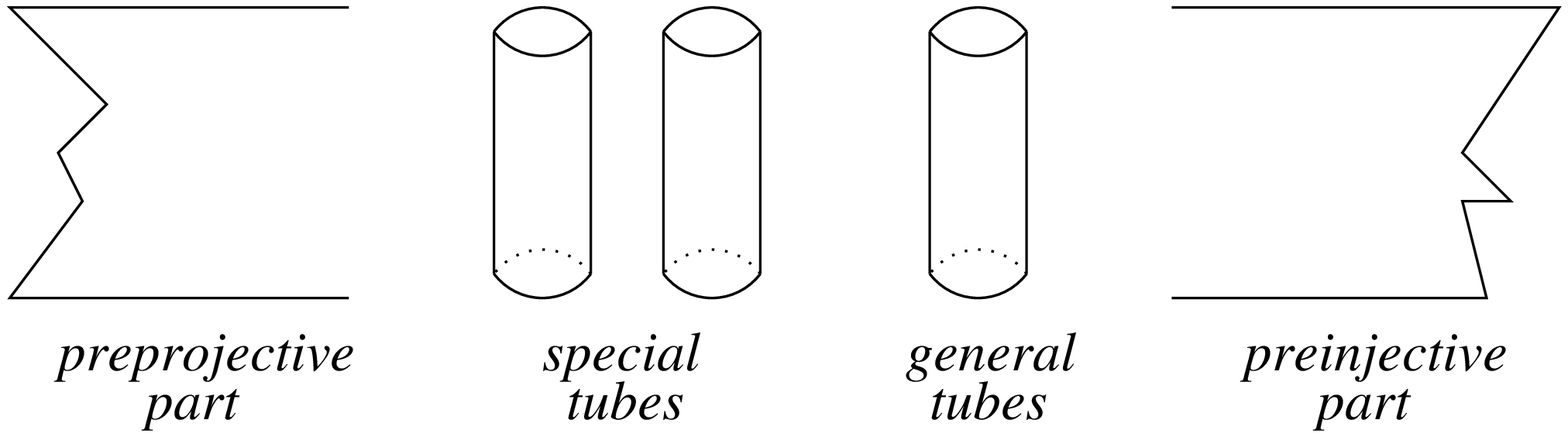}
\end{figure}

\clearpage
In this concrete case 

\begin{enumerate}

\item The preprojective series is:
\begin{figure}[ht]
\hspace{3cm}
\includegraphics[height=3cm,width=6cm,angle=0]{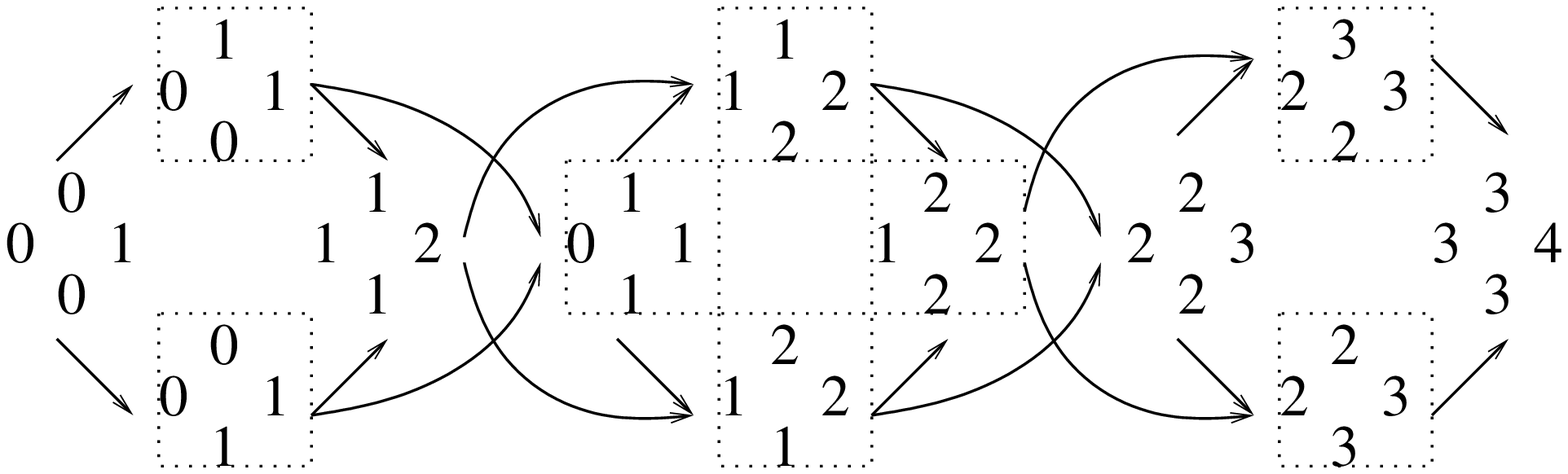}
\end{figure} 

(We mark with dotted boxes the objects that remain nonzero after tensoring by 
$\tilde{A}/I$.)

\item The preinjective series is
\begin{figure}[ht]
\hspace{3cm}
\includegraphics[height=3cm,width=6cm,angle=0]{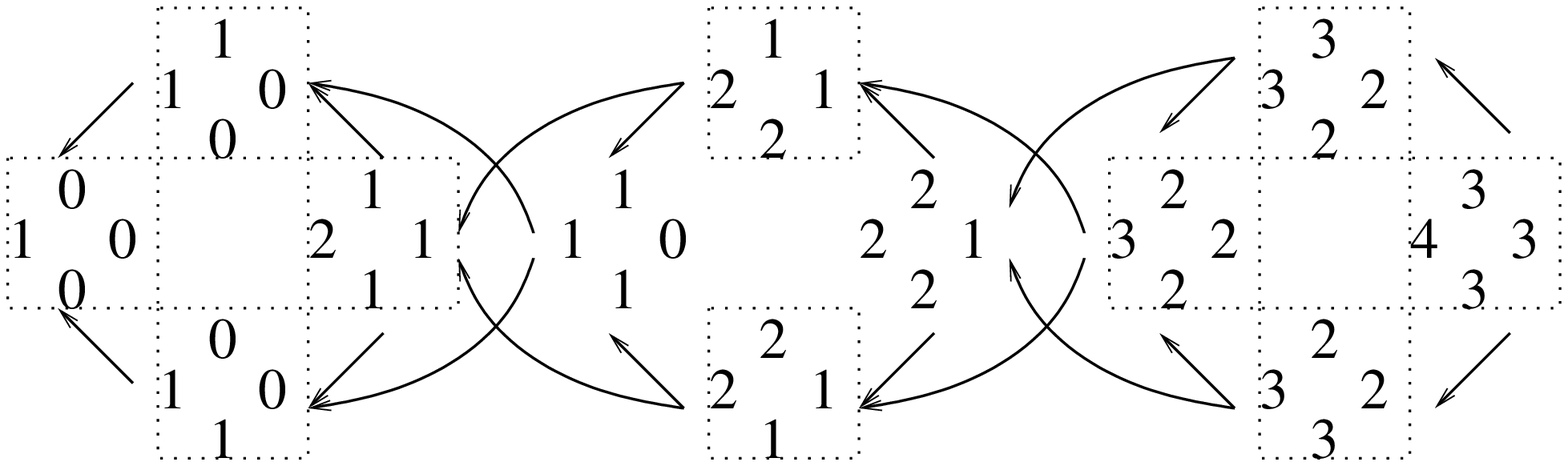}
\end{figure} 

\item Two special tubes are

\begin{figure}[ht]
\hspace{3cm}
\includegraphics[height=6cm,width=6cm,angle=0]{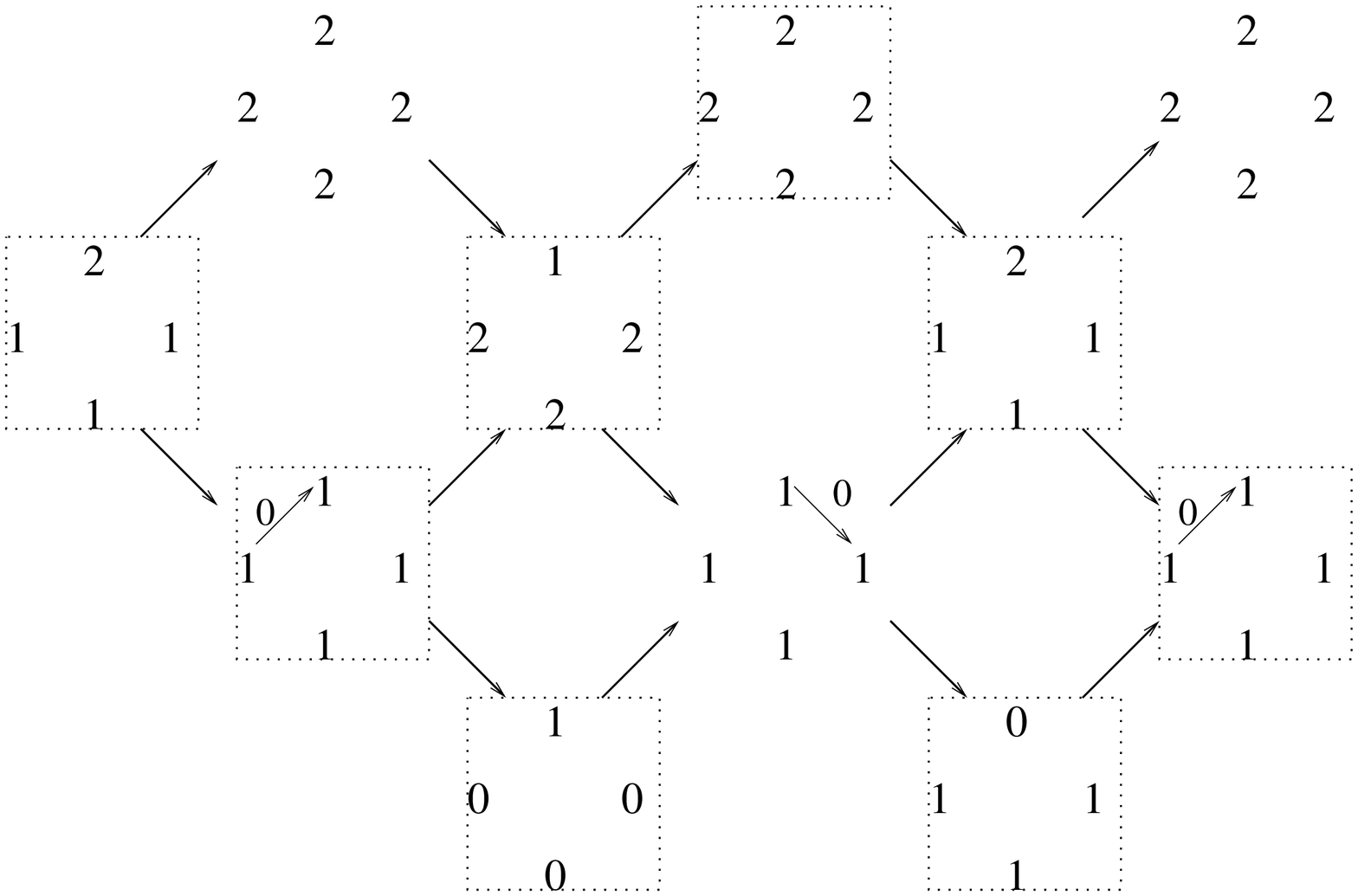}
\end{figure} 

and the symmetric one.
\end{enumerate}

We see that as preprojective series, as well as preinjective series and two special tubes are 2-periodic. Let $M$ be a preprojective module with the dimension vector $(d_1, d_2, d_3, d_4)$. Then $\tau^{-1}\circ \tau^{-1} (M)$ has the dimension vector $(d_1+2, d_2+2, d_3+2, d_4+2)$. The same holds for a preinjective module $N$ and  $\tau\circ \tau (N)$.  The same holds also for special tubes, if one goes two floors upstairs.

Consider the module $M$ from the preprojective series 

\begin{tabular}{p{3.7cm}c}
&
\xymatrix
{
& 2 \ar[rd] & \\
2 \ar[ru] \ar[rd] & & 3 \\
& 2 \ar[ru] & \\
}
\end{tabular} 

\noindent 
It has  a projective resolution 
$$
0 \lar P_4   \lar P_1^2 \lar M \lar 0.
$$
From it follows that   $\tilde{A}/I\tens_{\tilde{A}} M = 0$ . 
With other words, the module $M$ plays no role in our matrix problem. 
 
The list of modules, which are relevant for our problem is the following:

\begin{enumerate}
\item
Preprojective modules:

\begin{tabular}{p{0.5cm}ccc}
&
\xymatrix
{
& \kk^{n+1} \ar[rd] & \\
\kk^{n} \ar[ru] \ar[rd] & & \kk^{n+1} \\
& \kk^{n} \ar[ru] & \\
}
&
\medskip&
\xymatrix
{
& \kk^{n} \ar[rd] & \\
\kk^{n} \ar[ru] \ar[rd] & & \kk^{n+1} \\
& \kk^{n+1} \ar[ru] & \\
}
\end{tabular}

\noindent 
and 

\begin{tabular}{p{3.3cm}c}
&
\xymatrix
{
& \kk^{n+1} \ar[rd] & \\
\kk^{n} \ar[ru] \ar[rd] & & \kk^{n+1} \\
& \kk^{n+1} \ar[ru] & \\
}
\end{tabular}

\item
 Preinjective modules 

\begin{tabular}{p{0.5cm}ccc}
&
\xymatrix
{
& \kk^{n+1} \ar[rd] & \\
\kk^{n+1} \ar[ru] \ar[rd] & & \kk^{n} \\
& \kk^{n} \ar[ru] & \\
}
&
\medskip
&
\xymatrix
{
& \kk^{n} \ar[rd] & \\
\kk^{n+1} \ar[ru] \ar[rd] & & \kk^{n} \\
& \kk^{n+1} \ar[ru] & \\
}
\end{tabular}

\noindent
and

\begin{tabular}{p{3.3cm}c}
&
\xymatrix
{
& \kk^{n} \ar[rd] & \\
\kk^{n+1} \ar[ru] \ar[rd] & & \kk^{n} \\
& \kk^{n} \ar[ru] & \\
}
\end{tabular}

\item
 Modules from special tubes: 

\begin{tabular}{p{0.5cm}ccc}
&
\xymatrix
{
& \kk^{n+1} \ar[rd] & \\
\kk^{n} \ar[ru] \ar[rd] & & \kk^{n} \\
& \kk^{n} \ar[ru] & \\
}
&
\medskip
&
\xymatrix
{
& \kk^{n} \ar[rd] & \\
\kk^{n} \ar[ru]^{J_n(0)} \ar[rd] & & \kk^{n} \\
& \kk^{n} \ar[ru] & \\
}
\end{tabular}

\noindent
and 

\begin{tabular}{p{3.3cm}c}
&
\xymatrix
{
& \kk^{n} \ar[rd] & \\
\kk^{n+1} \ar[ru] \ar[rd] & & \kk^{n+1} \\
& \kk^{n+1} \ar[ru] & \\
}
\end{tabular}

\noindent
and symmetric ones.

\end{enumerate}

We have to compute the images of these objects after applying the 
left derived functor $\tilde{A}/I \tens_{\tilde{A}}.$ In order to do it we have to consider  the minimal projective resolutions of all these modules and then apply the functor  $\tilde{A}/I \tens_{\tilde{A}}.$

\begin{enumerate}
\item Preprojective series.  From the module 
$$
0 \lar P_4^{n}  \lar P_1^n \oplus P_2 \lar 0
$$
we get
$$
0 \lar 0 \lar \kk(2) \lar 0;
$$
from the module  
$$
0 \lar P_4^{n}  \lar P_1^n \oplus P_3 \lar 0
$$
will be 
$$
0 \lar 0 \lar \kk(3) \lar 0;
$$
from the module 
$$
0 \lar P_4^{n+1}  \lar P_1^n \oplus P_2 \oplus P_3 \lar 0
$$
we get 
$$
0 \lar 0 \lar \kk(2) \oplus \kk(3) \lar 0.
$$
\item Preinjective series. From the module 
$$
0 \lar P_2 \oplus P_3 \oplus P_4^{n}  \lar P_1^{n+1} \lar 0
$$
we get 
$$
0 \lar \kk(2) \oplus \kk(3) \lar 0 \lar 0;
$$
from the module 
$$
0 \lar P_2 \oplus  \ P_4^{n}  \lar P_1^{n+1} \lar 0
$$
will be 
$$
0 \lar \kk(2) \lar 0 \lar 0;
$$
from the module 
$$
0 \lar P_3  \oplus P_4^{n}  \lar P_1^{n+1} \lar 0
$$
will be 
$$
0 \lar \kk(3) \lar 0 \lar 0;
$$
\item Let us finally consider modules from special tubes.  
From the module 
$$
0 \lar P_4^{n+1}  \lar P_1^{n}\oplus P_2 \lar 0
$$
we get 
$$
0 \lar 0  \lar \kk(2) \lar 0;
$$
from the module 
$$
0 \lar P_4^{n}\oplus P_2  \lar P_1^{n}\oplus P_2 \lar 0
$$
will be 
$$
0 \lar \kk(2)  \stackrel{0}\lar \kk(2) \lar 0;
$$
from the module  
$$
0 \lar P_4^{n}\oplus P_2  \lar P_1^{n+1} \lar 0
$$
we get
$$
0 \lar \kk(2)  \lar 0 \lar 0.
$$
The case of the second special tube is completely symmetric. 

\end{enumerate}

What are induced morphisms between all these modules after applying the functor  $\tilde{A}/I\tens_{\tilde{A}}$? The image of the preprojective series is: 
 
\begin{tabular}{p{0.5cm}c}
&
\xymatrix
{
 \kk(2) \ar[rd] &    & \kk(3) \ar[rd] & &\kk(2) &\dots\\
  & \kk(2)\oplus \kk(3) \ar[ru] \ar[rd] & &\kk(2)\oplus \kk(3) \ar[ru] \ar[rd] 
& & \dots\\
\kk(3) \ar[ru] &  &\kk(2) \ar[ru] & & \kk(3) & \dots\\
}
\end{tabular}

We want to prove that all induced morphisms in this diagram are non-zero. 
It is well-known (see, for instance \cite{Ringel}) that all morphisms between preprojective modules are determined by the Auslander-Reiten quiver.  So, every morphism is a linear combination of 
finite paths of irreducible morphisms.  Consider the morphism

\begin{tabular}{p{1cm}ccc}
&
\xymatrix
{
& \kk \ar[rd] & &\\
0 \ar[ru] \ar[rd] & & \kk  \ar[r] &\\
& 0 \ar[ru] & & \\
}
&
\xymatrix
{
& \kk^{n+1} \ar[rd] & \\
\kk^n \ar[ru] \ar[rd] & & \kk^{n+1} \\
& \kk^n \ar[ru] & \\
}
\end{tabular}

It is clear that this map remains non-zero after applying 
 $\tilde{A}/I\tens_{\tilde{A}}.$ 
But from this fact follows that all morphisms we are interested in are non-zero 
after applying $\tilde{A}/I\tens_{\tilde{A}}.$  The same argument can be applied to preinjective modules. 

Consider finally the case of special tubes. 
The module

\begin{tabular}{p{3.3cm}c}
&
\xymatrix
{
& \kk^{n} \ar[rd] & \\
\kk^{n} \ar[ru]^{J_n(0)} \ar[rd] & & \kk^{n} \\
& \kk^{n} \ar[ru] & \\
}
\end{tabular}

\noindent
has a resolution
$$
0 \lar P_4^{n}\oplus P_2  \lar P_1^{n}\oplus P_2 \lar 0.
$$
It is indecomposable and its endomorphism algebra is local. Any endomorphism of this module is therefore either invertible or nilpotent. An isomorphism induces a map of the form:

\begin{tabular}{p{3.6cm}c}
&
\xymatrix
{
\kk(2) \ar[r]^{0} \ar[d]^{\lambda} & \kk(2) \ar[d]^{\lambda} \\
\kk(2) \ar[r]^{0} & \kk(2) \\
}
\end{tabular}

Let  $f$ be a nilpotent endomorphism. Consider an induced map of its projective resolution $f_\bullet$: 

\begin{tabular}{p{1.8cm}c}
&
\xymatrix
{0 \ar[r] & P_4^{n}\oplus P_2  \ar[r] \ar[d]^{f_1} & P_1^{n}\oplus P_2 
\ar[r] \ar[d]^{f_0}&  0 \\
0 \ar[r] & P_4^{n}\oplus P_2  \ar[r] & P_1^{n}\oplus P_2 \ar[r] &  0. \\
}
\end{tabular}

The map  $f_\bullet^n$ is homotopic to the zero map. Consider the component 
  $f_{1}^n |_{P_2} : P_2 \lar P_2$ of the map 
 $f_{1}^n$. It is zero modulo the radical and hence  is equal to zero itself. 
But then 
 $f_{1} |_{P_2} : P_2 \lar P_2$ is also zero. The same holds of coarse for 
 $f_{0} |_{P_2} : P_2 \lar P_2.$ 
We have shown that nilpotent morphisms induce the zero map  modulo 
  $\tilde{A}/I \tens_{\tilde{A}}$.

Finally observe that the chain of morphisms 

\begin{tabular}{p{2.9cm}c}
&
\xymatrix
{
& \kk^{n+1} \ar[rd] &  & \\
\kk^{n} \ar[ru] \ar[rd] & & \kk^{n} \\
& \kk^{n} \ar[ru] \ar[d]& & \\
& & & \\
}
\end{tabular}

\begin{tabular}{p{3.1cm}c}
&
\xymatrix
{
& \kk^{n} \ar[rd] & & \\
\kk^{n} \ar[ru]^{J_n(0)} \ar[rd] & & \kk^{n} \\
& \kk^{n} \ar[ru] \ar[d]& & \\
& & & \\
}
\end{tabular}

\begin{tabular}{p{2.6cm}c}
&
\xymatrix
{
& \kk^{n} \ar[rd] & \\
\kk^{n-1} \ar[ru] \ar[rd] & & \kk^{n-1} \\
& \kk^{n-1} \ar[ru] & \\
}
\end{tabular}

\noindent  induces modulo  $\tilde{A}/I\tens_{\tilde{A}}$ the following maps:

\begin{tabular}{p{3.6cm}c}
&
\xymatrix
{
 0 \ar[r] \ar[d]& \kk(2) \ar[d]^{1} \\
 \kk(2) \ar[r]^{0} \ar[d] & \kk(2) \ar[d]^{1} \\
 0 \ar[r]  & \kk(2). \\
}
\end{tabular}

\medskip
In the same way the chain of morphisms

\begin{tabular}{p{2.9cm}c}
&
\xymatrix
{
& \kk^{n-1} \ar[rd] &  & \\
\kk^{n} \ar[ru] \ar[rd] & & \kk^{n} \\
& \kk^{n} \ar[ru] \ar[d]& & \\
& & & \\
}
\end{tabular}

\begin{tabular}{p{3.1cm}c}
&
\xymatrix
{
& \kk^{n} \ar[rd] & & \\
\kk^{n} \ar[ru]^{J_n(0)} \ar[rd] & & \kk^{n} \\
& \kk^{n} \ar[ru] \ar[d]& & \\
& & & \\
}
\end{tabular}

\begin{tabular}{p{2.6cm}c}
&
\xymatrix
{
& \kk^{n} \ar[rd] & \\
\kk^{n+1} \ar[ru] \ar[rd] & & \kk^{n+1} \\
& \kk^{n+1} \ar[ru] & \\
}
\end{tabular}

\noindent
induce the maps

\begin{tabular}{p{3.6cm}c}
&
\xymatrix
{
 \kk(2) \ar[r] \ar[d]^{1} & 0 \ar[d] \\
 \kk(2) \ar[r]^0 \ar[d]^{1} & \kk(2) \ar[d] \\
 \kk(2) \ar[r]  & 0 \\
}
\end{tabular}

We have the same picture for the second symmetric tube. 
Now we observe that the matrix problem describing the derived
category $D^{-}(A-\mathop{\rm mod})$ is given by the following bunch of chains
(see \cite{mp}, \cite{CB}),
 where small circles correspond to the horizontal stripes, small rectangles below to the vertical stripes, dotted lines show the related stripes and arrows describe the ordering in the chains or, the same, possible transformations between
different horizontal stripes:

\begin{figure}[ht]
\hspace{0.8cm}
\includegraphics[height=13cm,width=10cm,angle=0]{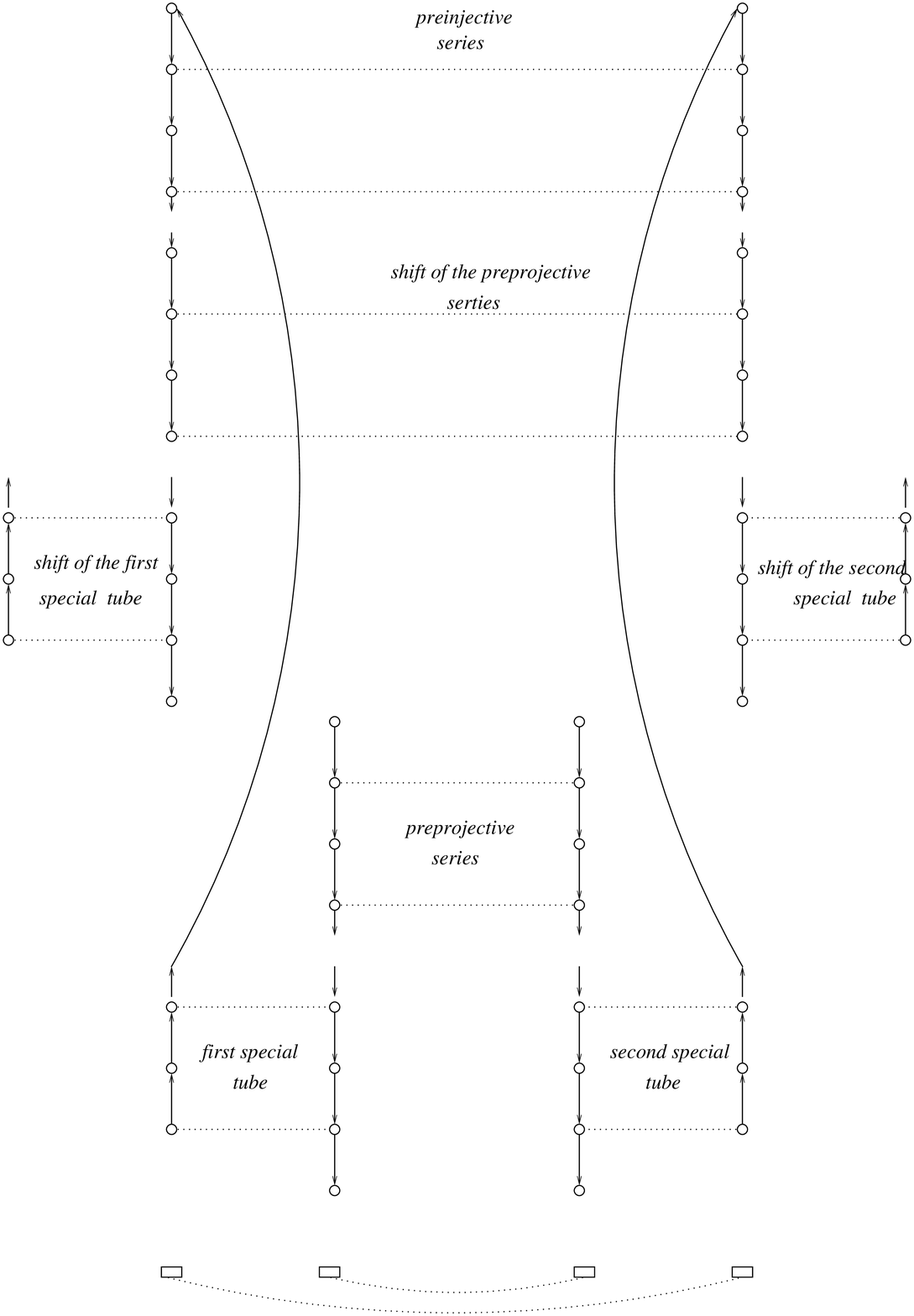}
\end{figure} 

\begin{enumerate}
\item We can do any simultaneous elementary transformations of  columns of 
the matrices
 $H_{k}(i)(0)$ and $H_{k}(i)(\infty)$.
\item We can do any simultaneous transformations of rows inside related blocks.
\item We can add a scalar multiple of any row from a block 
with lower weight to any row of a block of a higher weight (inside 
 the big matrix, of course). These transformations can be done independently
inside 
 $H_{k}(i)(0)$ and $H_{k}(i)(\infty)$.
\end{enumerate}

This type of matrix problems is well-known in the representation theory.
First they appeared in the work of Nazarova-Roiter (\cite{Roiter1}) about the 
classification of $\kk[[x,y]]/(xy)$-modules. They are called,
sometimes, Gelfand problems in honor of I.~M.~Gelfand, 
since they originated in a problem that first appeared in Gelfand's 
investigation of Harish-Chandra modules over $SL_{2}({\mathbb R})$ 
\cite{Gelfand} (see also \cite{mp}, \cite{CB}  and \cite{Jalgebra}).

\section{Derived categories of skew-gentle algebras}
In the same way as the representation theory of gentle algebras is based on 
 the
representation theory of hereditary algebras $\tilde{A}_{n}$, that of 
skew-gentle  algebras is
built on representations of $\tilde{D}_{n}$.

Consider the following example:

\begin{figure}[ht]
\hspace{4.3cm}
\includegraphics[height=2.2cm,width=3cm,angle=0]{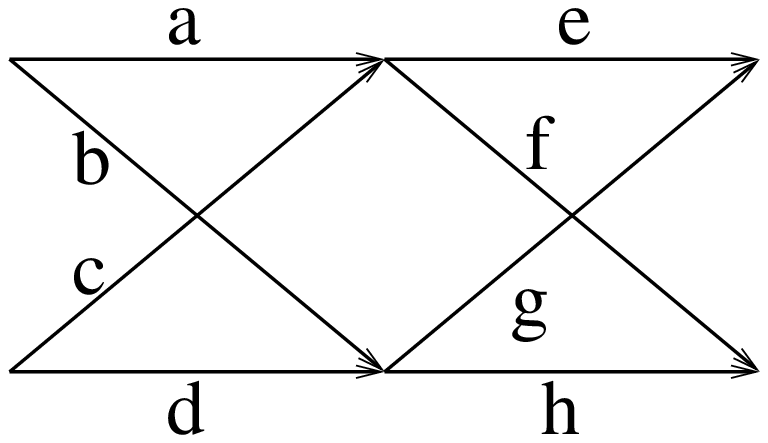}
\end{figure} 

\noindent
where all squares are commutative. Observe that we can embed this algebra
into 

\begin{figure}[ht]
\hspace{4.3cm}
\includegraphics[height=2.5cm,width=3cm,angle=0]{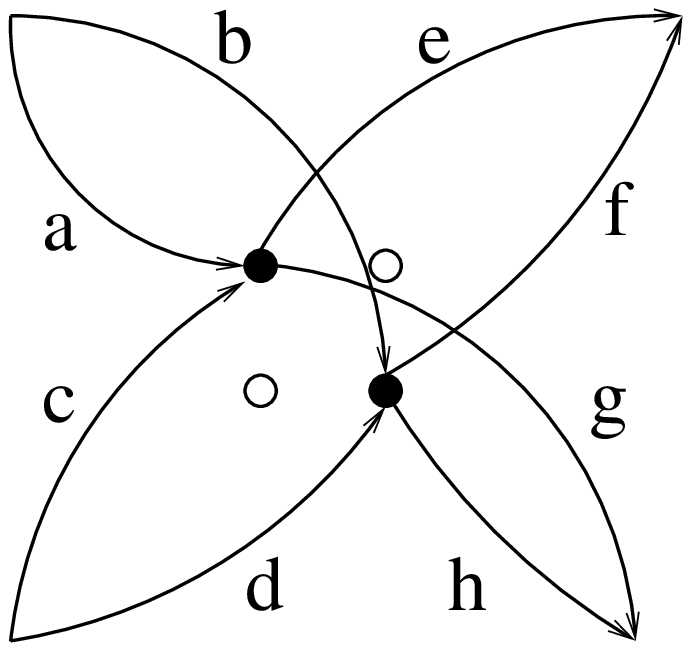}
\end{figure} 

\noindent
where $M_{2}(\kk)$ stands in the middle.  The last algebra is 
Morita-equivalent to the

\begin{figure}[ht]
\hspace{4.3cm}
\includegraphics[height=2.5cm,width=3cm,angle=0]{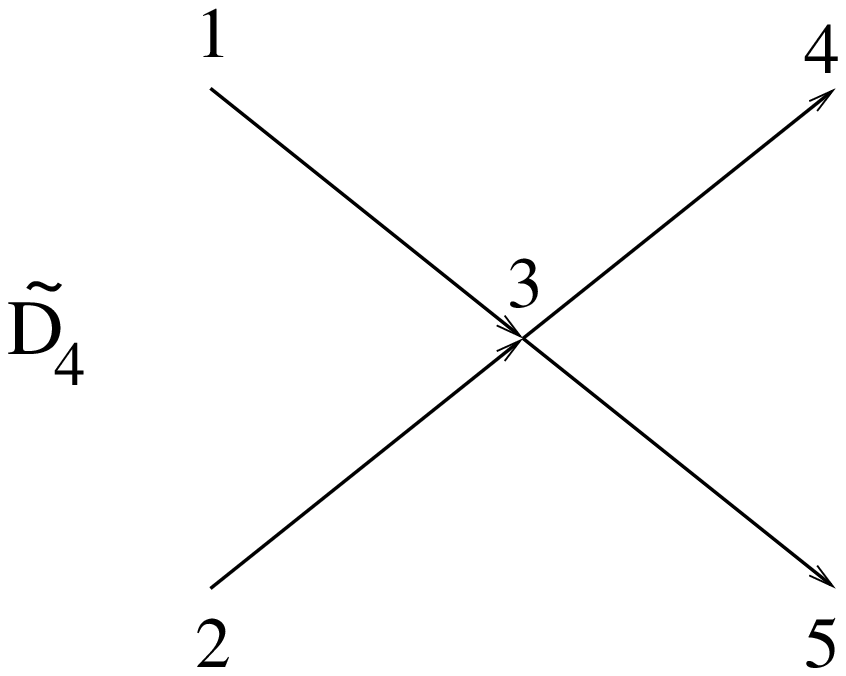}
\end{figure} 

\noindent
In this case the embedding $A/I \longrightarrow \tilde{A}/I$ is 
$\kk\times \kk \longrightarrow M_{2}(\kk)$.  It  means that the matrix problem we obtain is
the so called ``representation of the bunch of semi-chains'' (see \cite{mp},
\cite{CB} and \cite{Jalgebra} for more details).  
We have $\tilde{A}\otimes_{{\tilde A}/I} P_{i} = 0$, $i=1, 2, 4, 5.$

The continuous series of representations of $\tilde{D}_{4}$ has a projective
resolution 
$$
0 \longrightarrow P_{4}^n \oplus P_{5}^n \longrightarrow  P_{1}^n \oplus P_{2}^n \longrightarrow 
M_{n}(\lambda) \longrightarrow 0,
$$ 
hence $\tilde{A}/I\otimes_{\tilde A} M_{n}(\lambda) = 0.$
We again have to take care only of discrete series of $\tilde{D}_{4}$.
It consists of preprojective series, preinjective series, and three 
special tubes. 

\begin{enumerate}
\item
The preprojective series is:

\begin{figure}[ht]
\hspace{2.7cm}
\includegraphics[height=3cm,width=6cm,angle=0]{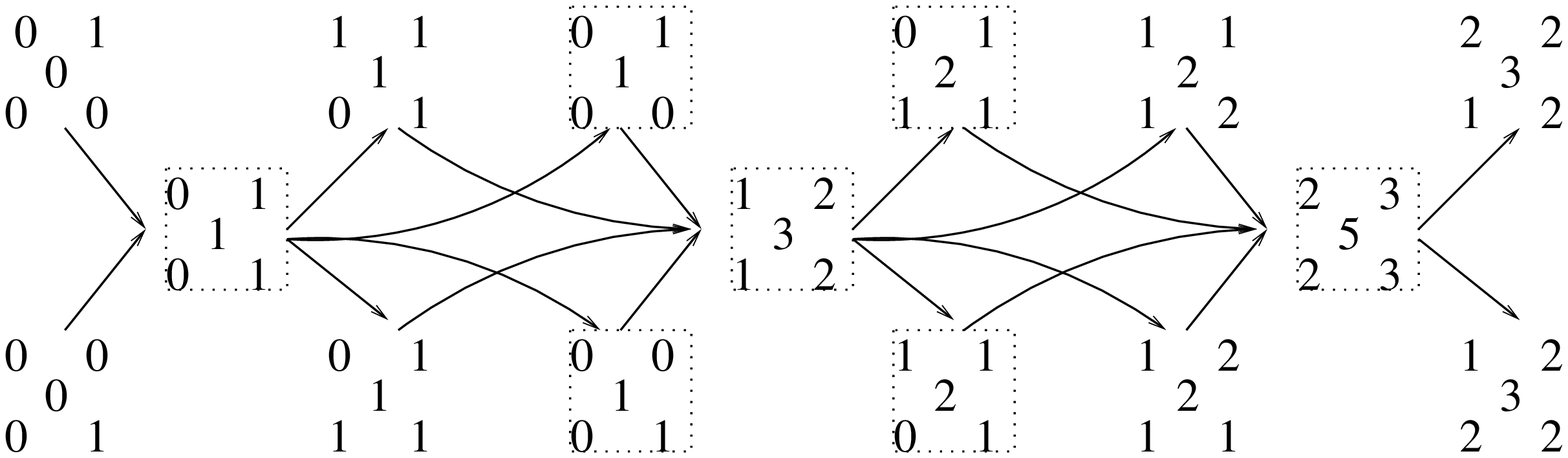}
\end{figure} 

\item The preinjective series is:

\begin{figure}[ht]
\hspace{2.7cm}
\includegraphics[height=3cm,width=6cm,angle=0]{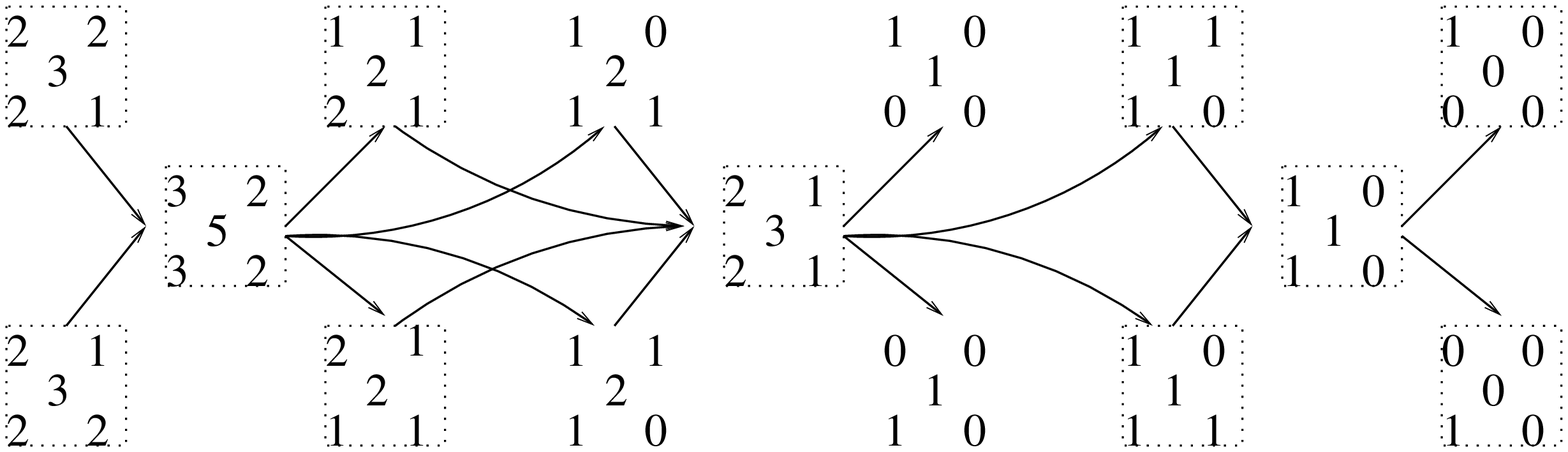}
\end{figure}

\item The only special tube that gives an input into the matrix problem is:

\begin{figure}[ht]
\hspace{3.8cm}
\includegraphics[height=4.5cm,width=4.5cm,angle=0]{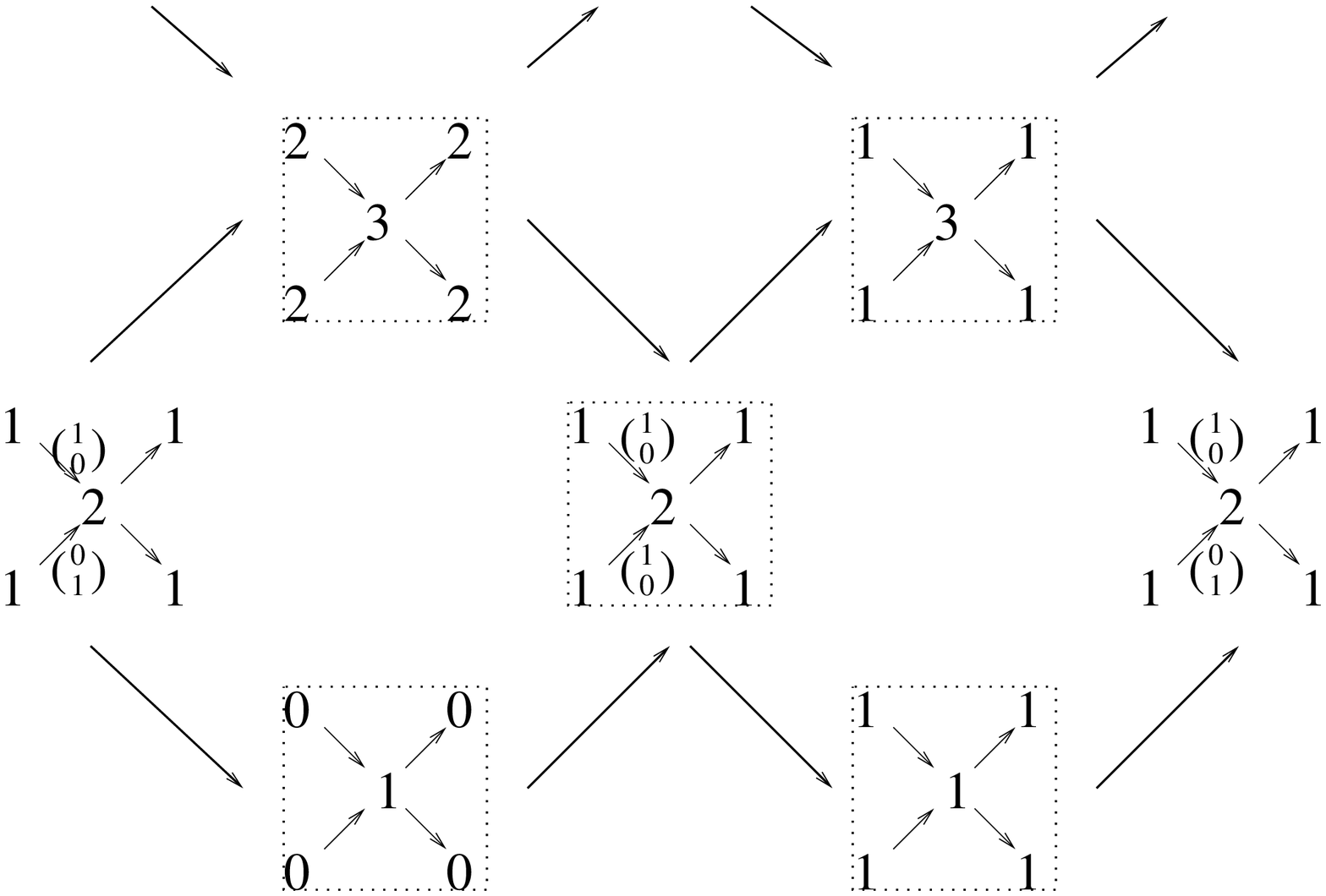}
\end{figure} 

\end{enumerate}
The representations from two other special tubes vanish under tensoring by 
$\tilde{A}/I$.

In the same way as for gentle algebras, we see that our matrix problem 
is given by the following partially ordered set:

\clearpage
\begin{figure}[ht]
\hspace{2cm}
\includegraphics[height=13cm,width=10cm,angle=0]{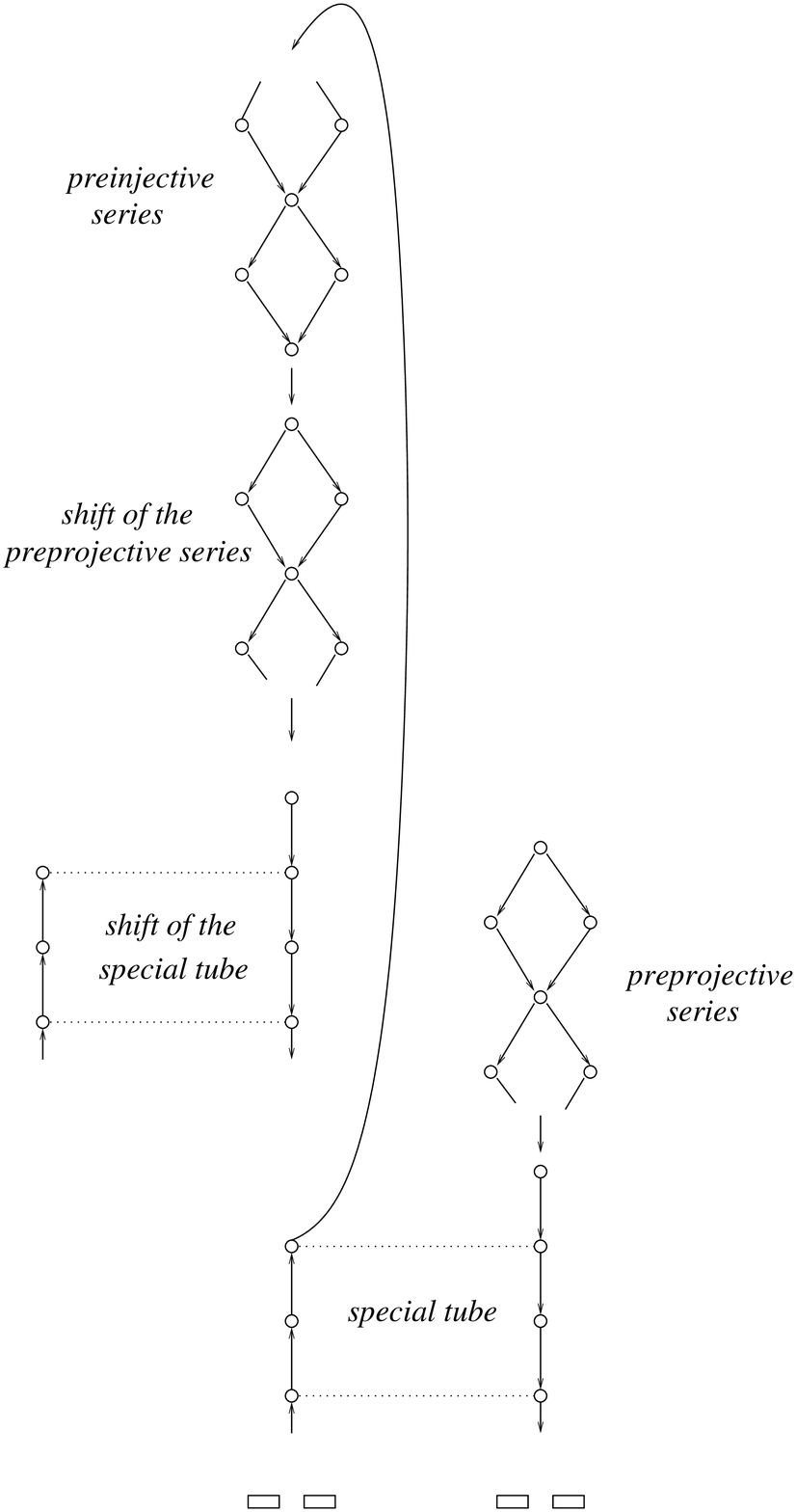}
\end{figure}

The generalization of this approach to other skew-gentle algebras of finite homological dimension gives a new proof of the result of 
Ch.~Gei\ss{} and J.~A.~de~la~Pena \cite{laPena} that the bounded derived 
categories of these  algebras are  tame.

\section{Skew-gentle algebras of infinite homological dimension}
Consider now our next example:

\clearpage

\begin{figure}[ht]
\hspace{4cm}
\includegraphics[height=1.5cm,width=4cm,angle=0]{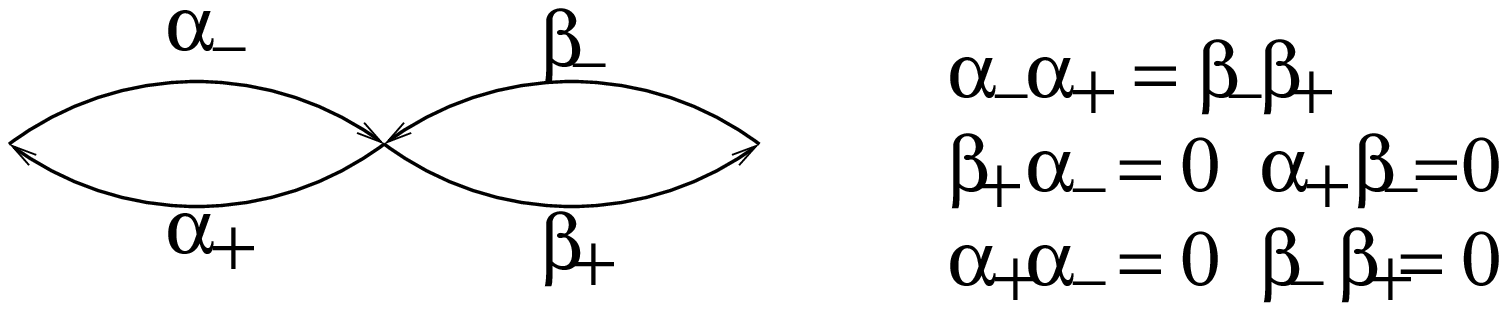}
\end{figure} 

\noindent
This algebra is skew-gentle of infinite homological dimension.
We can embed it in

\begin{figure}[ht]
\hspace{5cm}
\includegraphics[height=2cm,width=4cm,angle=0]{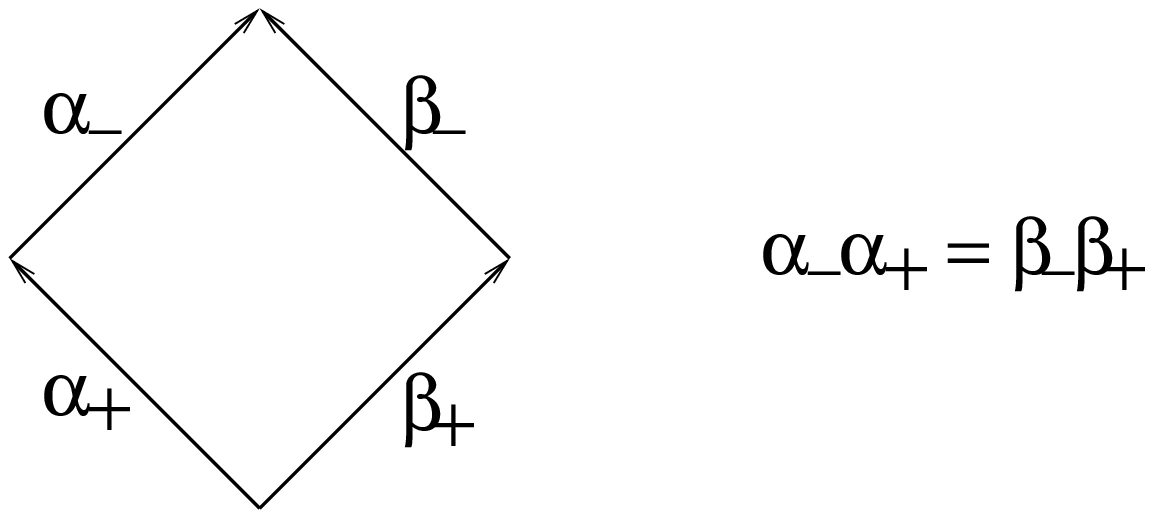}
\end{figure} 

It is well known (see, for example, \cite{Happel}) that this algebra is 
derived-equivalent to the algebra $\tilde{D}_{4}$.  In such a way we can
obtain our matrix problem. But we can embed it further into

\begin{figure}[ht]
\hspace{4.5cm}
\includegraphics[height=1cm,width=3cm,angle=0]{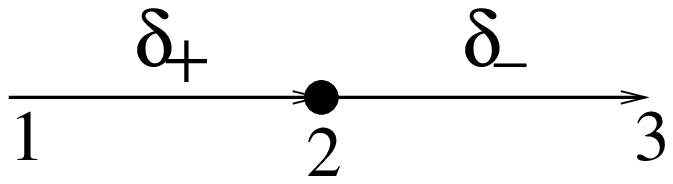}
\end{figure} 

\noindent
where the fat point in the middle means $M_{2}(\kk)$. This algebra is 
Morita-equivalent to the algebra $A_{2}$. Note that there are only the 
following indecomposable complexes in $D^{-}(A_{2}-\mathop{\rm mod})$ (up to shifts):
$$\dots 0\longrightarrow P_{1} \longrightarrow 0 \longrightarrow \dots,$$ 
$$\dots 0\longrightarrow P_{2} \longrightarrow 0 \longrightarrow \dots,$$
$$\dots 0\longrightarrow P_{3} \longrightarrow 0 \longrightarrow \dots;$$
$$\dots 0\longrightarrow P_{3} \longrightarrow P_{1} \longrightarrow 0 \longrightarrow \dots,$$
$$\dots 0\longrightarrow P_{3} \longrightarrow P_{2} \longrightarrow 0 \longrightarrow \dots,$$
$$\dots 0\longrightarrow P_{2} \longrightarrow P_{1} \longrightarrow 0 \longrightarrow \dots.$$

It is easy to establish a matrix problem now:

\begin{figure}[ht]
\hspace{2.5cm}
\includegraphics[height=3cm,width=8cm,angle=0]{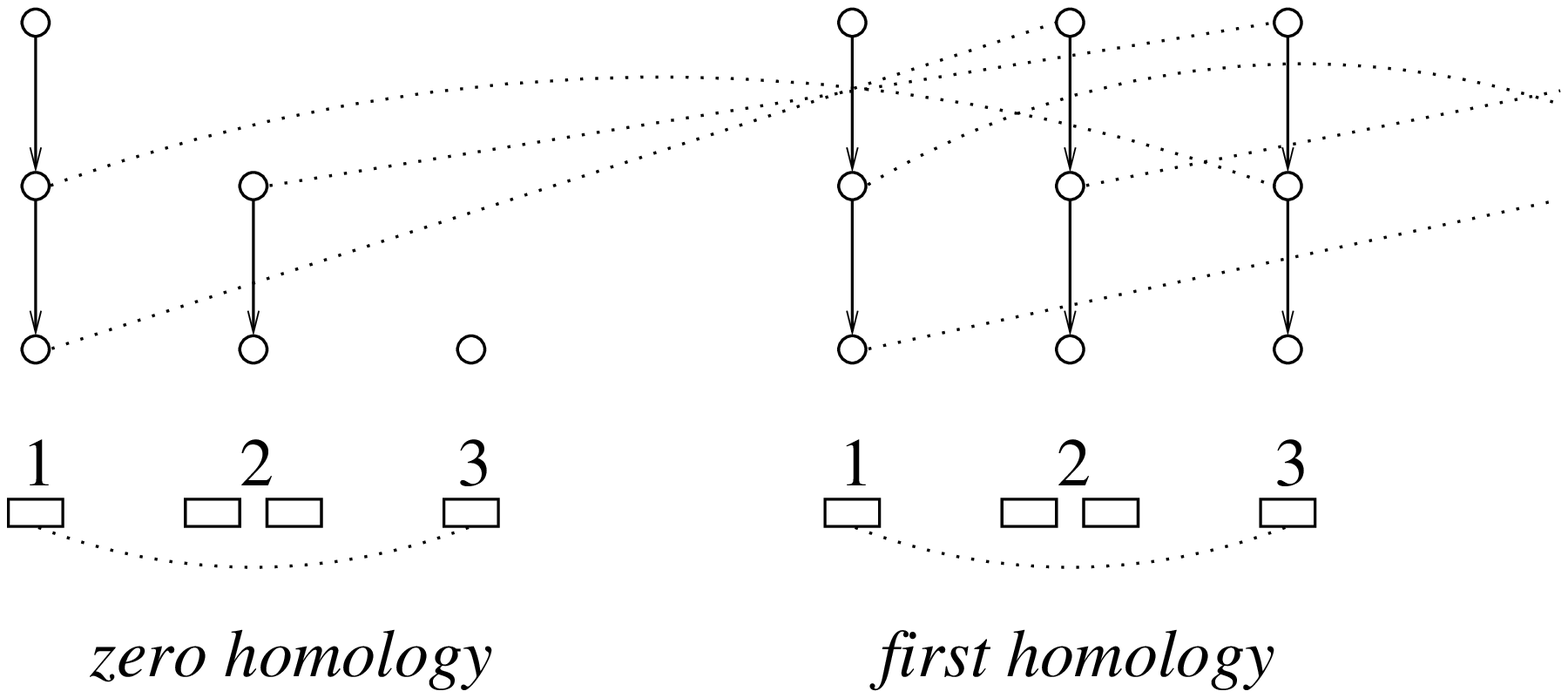}
\end{figure} 

Derived categories of skew-gentle algebras of infinite homological dimension were independently considered in \cite{Bekkert2}.

\section{Derived category of degenerated tubular algebra $(2,2,2,2; 0)$}

In all previous examples we embedded our path algebra $A$ into a hereditary algebra $\tilde{A}$.  As we shall see in the following example, it is also possible to consider an embedding into a concealed algebra.

If we consider a canonical algebra of tubular type $(2,2,2,2; \lambda)$ 
and set the forbidden value of parameter $\lambda = 0$, then we get the following quiver

\begin{figure}[ht]
\hspace{3.9cm}
\includegraphics[height=3cm,width=4cm,angle=0]{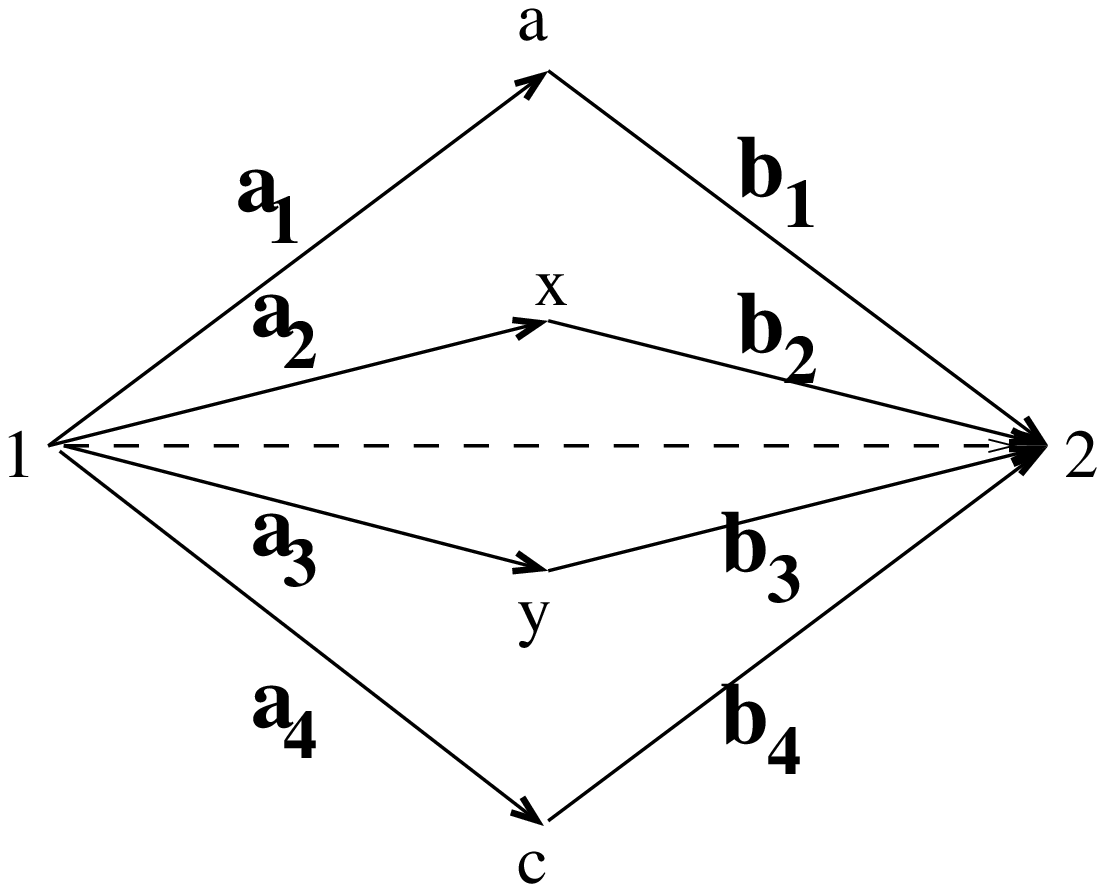}
\end{figure}

\noindent
with relations
$$
 b_2 a_2 = b_3 a_3
$$
and
$$
 b_1 a_1 + b_2 a_2 + b_4 a_4 = 0.
$$
We can do our trick with gluing of idempotents $x$ and $y$ and embed this algebra into

\begin{figure}[ht]
\hspace{3.9cm}
\includegraphics[height=2.5cm,width=4cm,angle=0]{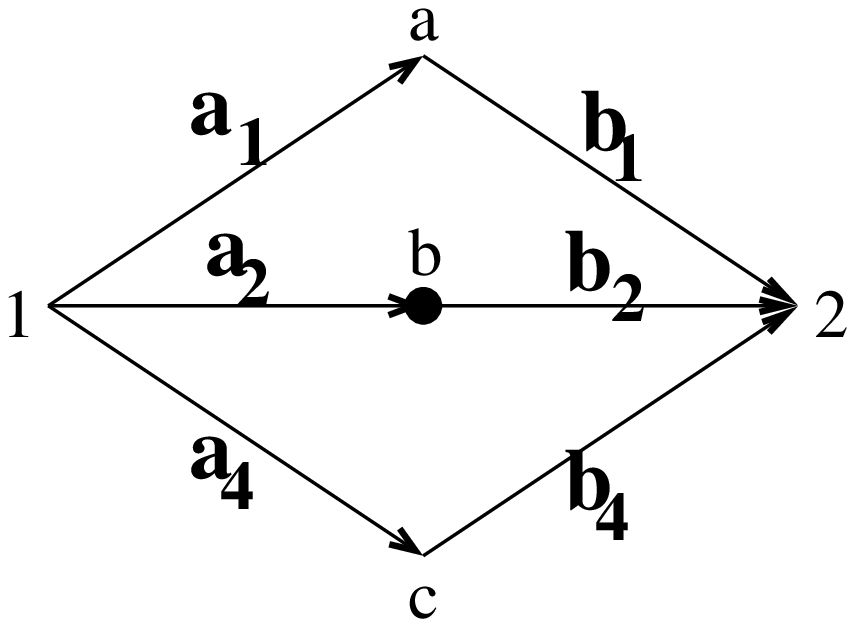}
\end{figure}

\noindent
where the fat point as usually means $M_2(\kk)$.   The corresponding basic algebra is a tame-concealed algebra of type $(2,2,2)$. It is well-known
\cite{Ringel} that it is derived-equivalent to  the algebra $\tilde{D}_4.$
Take the ideal $I$ equal to the ideal generated by all idempotents
$e_1, e_2, e_{a}, e_{c}$.  Then we have: $A/I = \kk\times \kk$, 
$\tilde{A}/I = Mat_2(\kk)$ and the map $A/I \lar \tilde{A}/I$ is the diagonal embedding.

It is no longer true that  any complex of $D^{b}(\tilde{A}-\mod)$ 
is isomorphic to its homology.  The structure of the  Auslander-Reiten quiver is the same as for derived categories of tame hereditary algebras.

\clearpage

\begin{figure}[ht]
\hspace{1cm}
\includegraphics[height=3cm,width=11cm,angle=0]{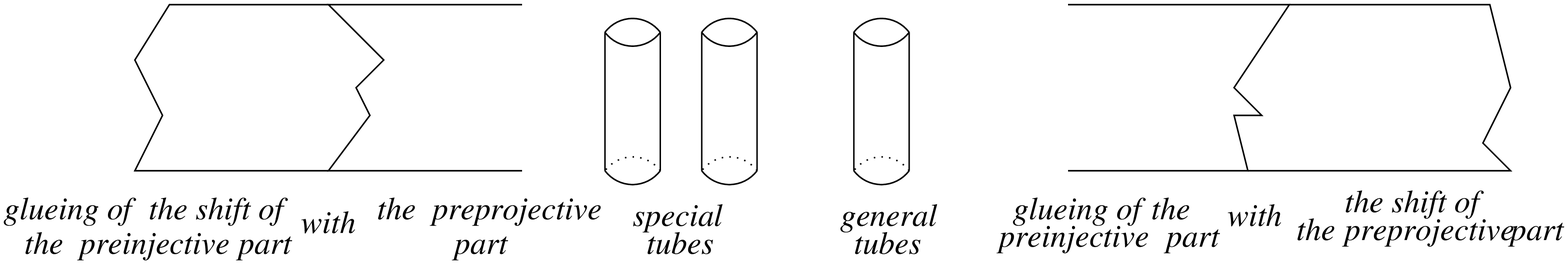}
\end{figure} 

In particular, continuous series of complexes are just shifts of modules  of tubular type. 
But  they have the following form:

\begin{tabular}{p{2.2cm}c}
&
\xymatrix
{
   && \kk^n \ar[rrd]^{I_n}&&  \\
\kk^n \ar[rr]^{I_n} \ar[rru]^{I_n} \ar[rrd]_{I_n} && \kk^n \ar[rr]_{-I_n - J_n(\lambda)}&& \kk^n \\
 && \kk^n \ar[rru]_{J_n(\lambda)}&&  \\
}
\end{tabular}

It is easy to see that they they minimal projective resolutions have the form
$ P_2^n \lar P_2^n$ and hence $\tilde{A}/I \tens_{\tilde{A}} 
(P_2^n \lar P_2^n) = 0$, they do not affect the resulting matrix problem.

Let us first consider the structure of preprojective and preinjective 
components of the Auslander-Reiten quiver of $\tilde{A}-\mod$.

\begin{figure}[ht]
\hspace{1.7cm}
\includegraphics[height=4cm,width=9cm,angle=0]{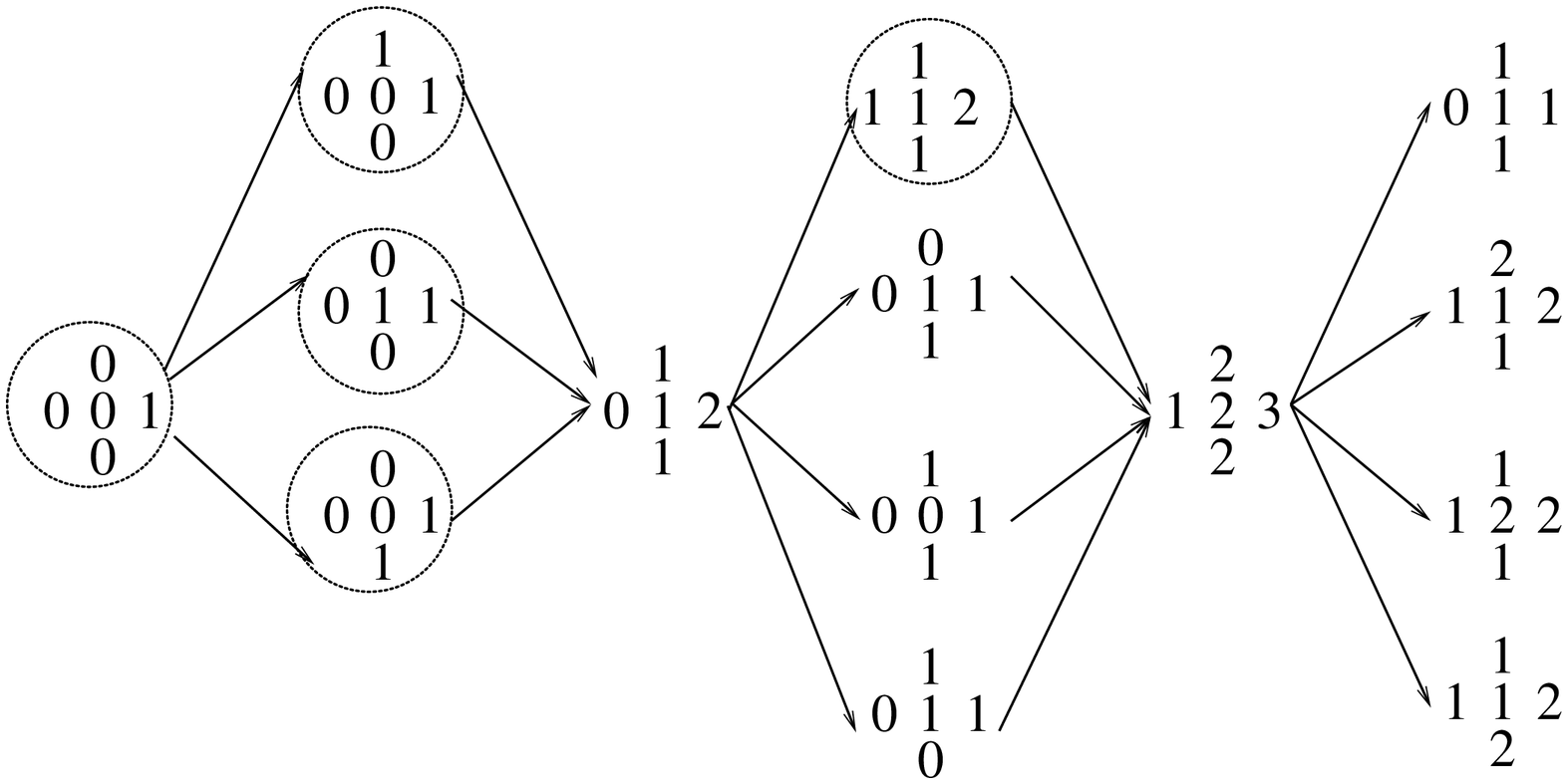}
\end{figure}

\begin{figure}[ht]
\hspace{1.7cm}
\includegraphics[height=4cm,width=9cm,angle=0]{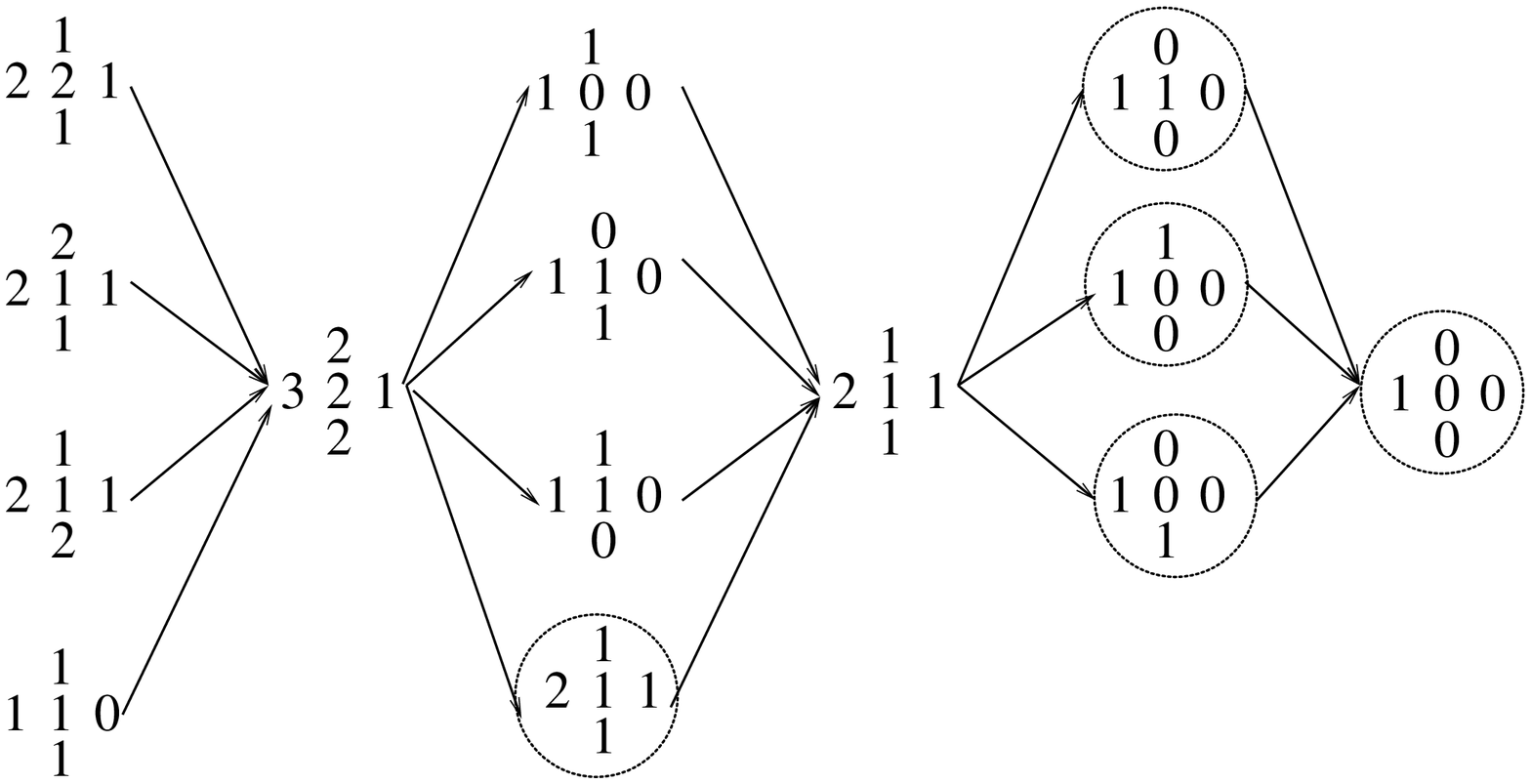}
\end{figure}

We can use the following lemma  (see \cite{Happel})

\begin{lemma}
Let $A$ be an associative $\kk$-algebra, 
$$
0 \lar M \stackrel{u}\lar N \stackrel{v}\lar K \lar 0
$$
an Auslander-Reiten sequence in $A-\mod$, $w \in \mbox{\rm Ext}^1(K, M) = 
\mbox{\rm Hom}_{D^{b}(A)}(K, T(M))$ the 
corresponding element. Then the following is equivalent:

\begin{enumerate}
\item 
$
M \stackrel{u}\lar N \stackrel{v}\lar  K \stackrel{w}\lar T(M)$ 
is an Auslander-Reiten triangle in \\ $D^{b}(A-\mod)$.
\item 
$
\mbox{\rm inj.dim}(M) \leq 1,$
$
\mbox{\rm proj.dim}(K) \leq 1.$
\item 
$
\mbox{\rm Hom}_A (I, M) = 0 
$
for any injective $A$-module $I$ and 
$
\mbox{\rm Hom}_A (K, P) = 0 
$
for any projective $A$-module $P$.
\end{enumerate}
\end{lemma}

This lemma means that the structure of the Auslander-Reiten quiver of the category $D^{b}(A-\mod)$ is basically the same as for $A-\mod$. For instance, 
all morphisms from tubes are still almost split in the derived category. 

There is exactly one indecomposable complex in $ D^{b}(A-\mod)$, which is not isomorphic to a shift of some module: it is 
$$
P_a \oplus P_b \oplus P_c \lar P_1. \footnote{
The first author would like to thank C.-M.Ringel for explaining him this.
}
$$ 
It is easy to see that this complex has two non-trivial homologies and is indecomposable. Now from lemma above and the fact that there is only one non-trivial indecomposable complex we can derived the exact form of the ``gluing'' of the pre-projective component with the shift of the pre-injective component in the 
Auslander-Reiten quiver:

\begin{figure}[ht]
\hspace{-1cm}
\includegraphics[height=4.2cm,width=14cm,angle=0]{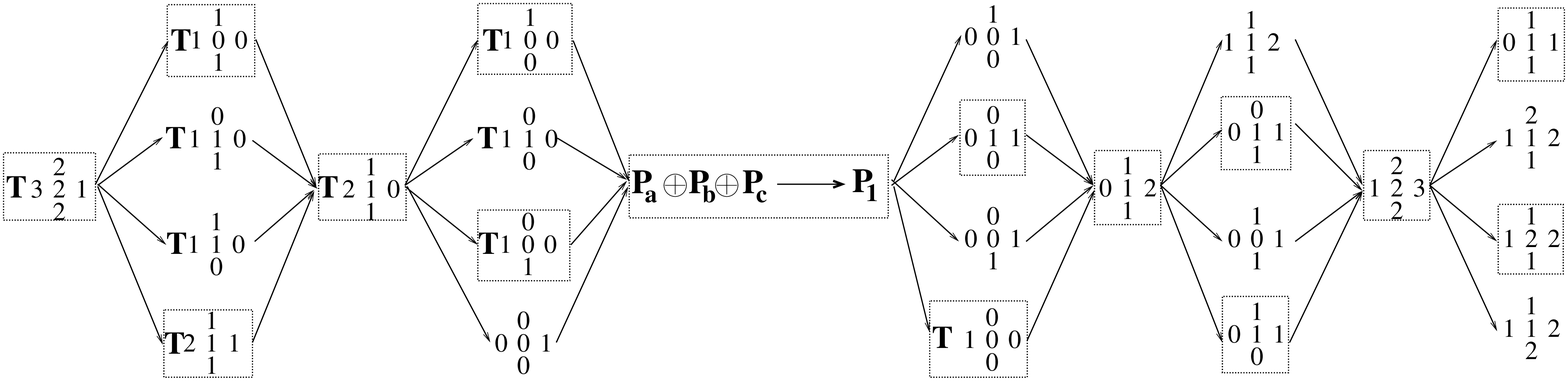}
\end{figure}

\clearpage

There are finally 3 special tubes of the length 2. They are completely symmetric
with respect to permutation of vertices $a$, $b$ and $c$ and only one of them is relevant for the matrix problem.

\begin{figure}[ht]
\hspace{3cm}
\includegraphics[height=6cm,width=6cm,angle=0]{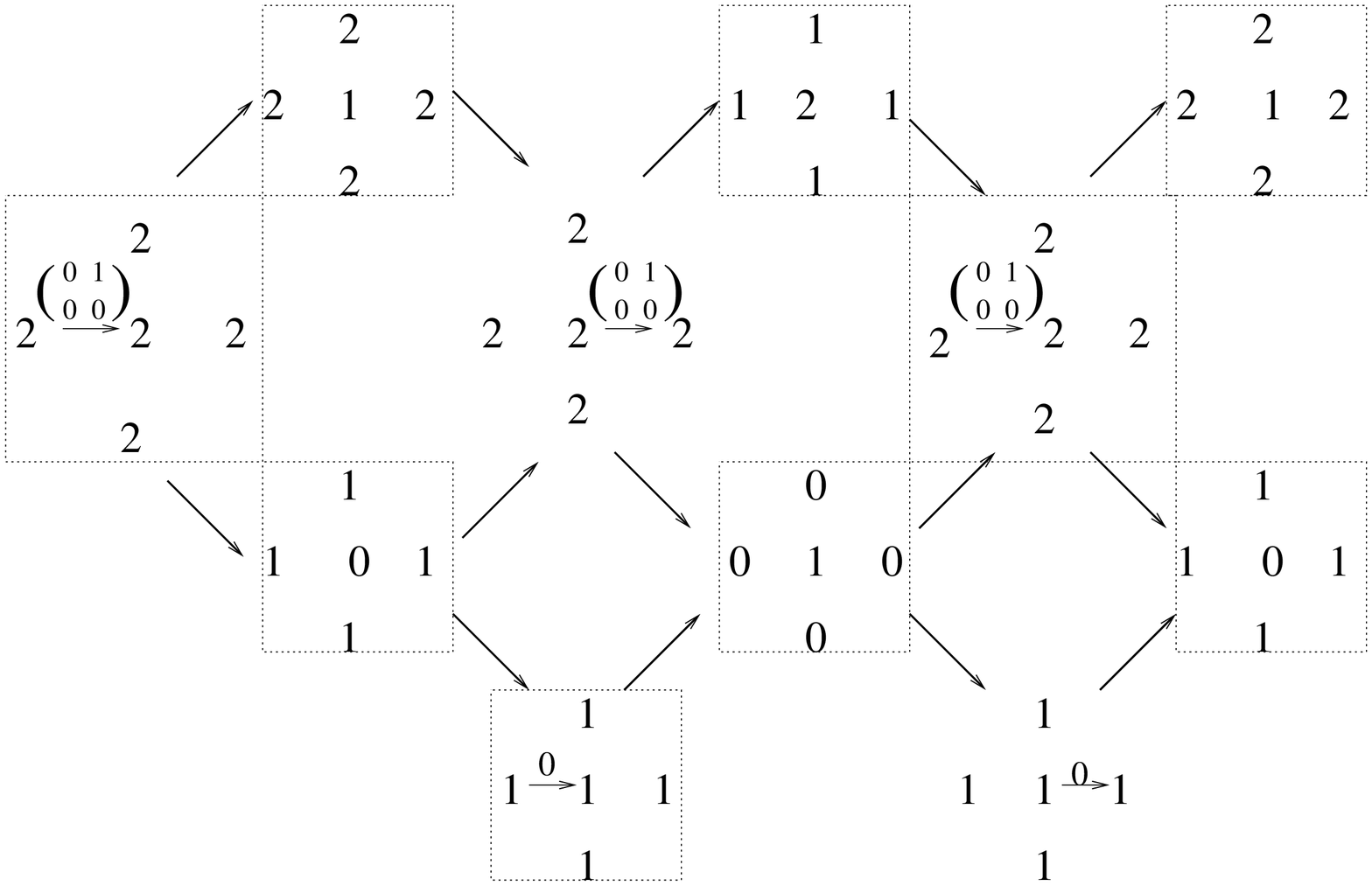}
\end{figure}

The resulting matrix problem is given by the following partially ordered set
(which is very similar to the case of skew-gentle algebra, considered above)

\begin{figure}[ht]
\hspace{2cm}
\includegraphics[height=9.5cm,width=6cm,angle=0]{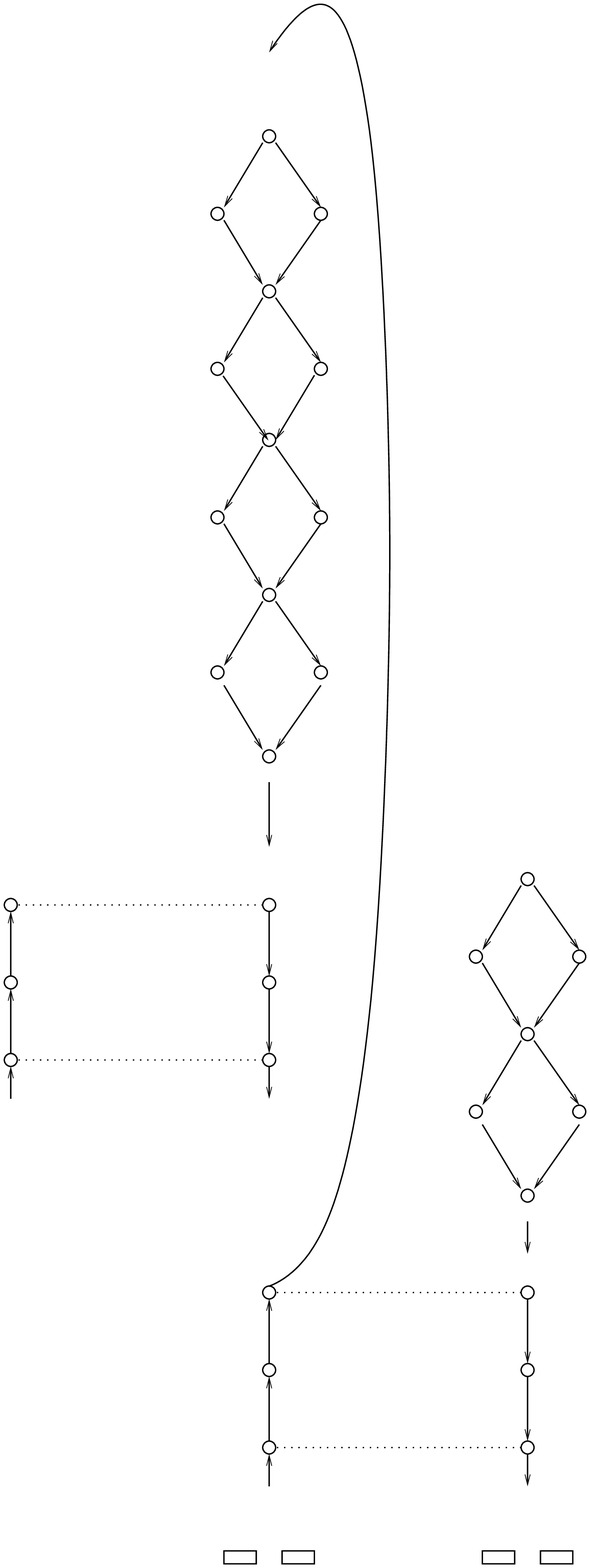}
\end{figure}

\clearpage
Therefore the degenerated tubular algebra $(2,2,2,2;0)$ is derived-tame of 
exponential growth. 

Is seems to be very plausible that this algebra is closely related to 
``weighted projective lines with a singularity of virtual genus one''
 and the map $A\lar \tilde{A}$
plays the role of non-commutative normalization (and note that $\tilde{A}$ correspond to a weighted projective line of virtual genus zero).  We are planing to come back to this problem in the future. 

One  can also consider a mixes situation where we consider a radical embedding of a $\kk$-algebra $A$ into a product of finite or tame hereditary algebras and
finite or tame concealed algebras.

%\vspace{1cm}
%\begin{tabular}{lp{1.5cm}l}
%Igor Burban & & Yuriy Drozd\\
%Universit\"at Kaiserslautern & & Kyiv University \\
%burban@mathematik.uni-kl.de & & yuriy@drozd.org
%\end{tabular}

\end{document}